\documentclass[a4paper]{article}
\usepackage{amsmath,amsfonts,amssymb,amsthm}
\usepackage{graphics}
\usepackage{epsfig}
\usepackage{color}
\usepackage[hang,flushmargin]{footmisc} 

\def\R{{\mathbb{R}}}
\def\Z{{\mathbb{Z}}}

\def\E{{\mathbb{E}}}

\newcommand{\overbar}[1]{\mkern 1.5mu\overline{\mkern-1.5mu#1\mkern-1.5mu}\mkern 1.5mu}

\DeclareMathOperator*{\osc}{osc}

\newcommand{\expec}[1]{\langle #1 \rangle}

\newcommand{\ignore}[1]{}

\newtheorem{definition}{Definition}
\newtheorem{proposition}{Proposition}
\newtheorem{theorem}{Theorem}
\newtheorem{remark}{Remark}
\newtheorem{lemma}{Lemma}
\newtheorem{corollary}{Corollary}
\newtheorem*{corollary*}{Corollary}

\title{Annealed estimates on the Green function}
\author{Daniel Marahrens\thanks{Max Planck Institute for Mathematics in the Sciences, Inselstr.\ 22, 04103 Leipzig,
Germany, \texttt{Daniel.Marahrens@mis.mpg.de} resp.\ \texttt{Felix.Otto@mis.mpg.de}} \and Felix Otto\footnotemark[1]}
\date{}

\begin{document}

\maketitle

\begin{abstract}
We consider a random, uniformly elliptic coefficient field $a(x)$ on the 
$d$-dimensional integer lattice $\mathbb{Z}^d$. 
We are interested in the spatial decay of the quenched elliptic Green function $G(a;x,y)$. 
Next to stationarity, we assume that the spatial correlation of the coefficient field
decays sufficiently fast to the
effect that a logarithmic Sobolev inequality holds for the ensemble $\langle\cdot\rangle$. 
We prove that
{\it all} stochastic moments of the first and second mixed derivatives of the Green
function, that is, $\langle|\nabla_x G(x,y)|^p\rangle$ and 
$\langle|\nabla_x\nabla_y G(x,y)|^p\rangle$,  
have the same decay rates in $|x-y|\gg 1$ as for the constant coefficient Green function,
respectively.
This result relies on and substantially extends the one by Delmotte and Deuschel \cite{DeuschelDelmotte}, which optimally controls
second moments for the first derivatives and first moments of the second mixed derivatives
of $G$, that is, $\langle|\nabla_x G(x,y)|^2\rangle$ and  
$\langle|\nabla_x\nabla_y G(x,y)|\rangle$. 
As an application, we are able to obtain optimal estimates on the random part of the homogenization error even for large
ellipticity contrast.
\end{abstract}

\renewcommand{\thefootnote}{\fnsymbol{footnote}} 
\footnotetext{\noindent\today\\\emph{MSC2010 subject classifications.} Primary: 35B27, Secondary: 35J08, 39A70, 60H25\\
	      \emph{Key words} Stochastic homogenization, elliptic equations, Green function, annealed estimates}     
\renewcommand{\thefootnote}{\arabic{footnote}} 

%


\section*{Outline}

The outline of this work is as follows: After introducing the discrete setting in Section \ref{sec:not}, we present the statistical
assumptions and the main result on the annealed moments of the Green function in Section \ref{sec:main}. The following two sections
contain applications of the main result: We present optimal estimates on the random part of the homogenization error in
Section~\ref{sec:hom} and Section \ref{sec:DeGiorgi} contains an annealed H\"older-estimate in the spirit of De Giorgi. In Section~\ref{sec:LSI} we explain our main assumption, a logarithmic Sobolev inequality,
which in particular holds for all independent, identically distributed coefficient fields. Section \ref{sec:ingredients} contains the main
ingredients of the proof of the annealed Green function estimates --- in particular we recall the result by Delmotte and Deuschel
\cite{DeuschelDelmotte}. All proofs are postponed until Section \ref{sec:proofs}.

\section{Discrete uniformly elliptic equations}\label{sec:not}

In this paper we consider linear second-order difference equations with uniformly elliptic, bounded random coefficients of the form
\begin{equation}\label{PDE}
 \nabla^*(a \nabla u)(x) = f(x) \quad \text{for all $x\in\Z^d$.}
\end{equation}
If there is no danger of confusion, we also write $\nabla^* a \nabla u$ for $\nabla^*(a\nabla u)$. 
In this equation we define the \textit{spatial derivatives} as follows: Let $\E^d$ denote the set of \emph{edges} of $\Z^d$ consisting of
all pairs $[x,x+e_i]$ of neighboring vertices with $x\in\Z^d$, $i=1,\ldots,d$, where $e_1,\ldots,e_d$ is the canonical basis of $\R^d$.
For functions on vertices $\zeta:\Z^d\to\R$ and functions on edges 
$\xi:\E^d\to\R$ we set
\begin{align*}
 \nabla \zeta ([x,x+e_i]) &= \zeta(x+e_i)-\zeta(x),\\
 \nabla^*\xi(x) &= \sum_{i=1}^d \big(\xi([x-e_i,x]) - \xi([x,x+e_i])\big).
\end{align*}
The spatial derivatives $\nabla \zeta$ and $-\nabla^*\xi$ are the discrete gradient and divergence, respectively, on the lattice $\Z^d$.
As our notation suggests, the operators $\nabla$ and $\nabla^*$ are adjoint in the sense of
\[
 \sum_{e\in\E^d} \xi(e) \nabla\zeta(e) = \sum_{x\in\Z^d} \nabla^*\xi(x) \zeta(x).
\]
In \eqref{PDE}, the  coefficient field  $a$ is a field on edges $a:\E^d \to \R$. Consequently $\nabla^* a \nabla$ is well-defined as an operator on vertex fields $\Z^d\to\R$. In this paper, we denote edges in $\E^d$ by the
letters $e$ and $b$ and vertices in $\Z^d$ by the letters $x$, $y$, and $z$.

\medskip

Throughout this work we consider coefficient fields $a:\E^d \to \R$ in the space $\Omega$ of \emph{uniformly elliptic} coefficient fields, i.e.\ we let
\begin{equation}\label{Omega}
 \Omega:= \big\{ a:\E^d \to \R : \lambda \le a(e) \le 1 \text{ for all $e\in\E^d$} \big\} = [\lambda,1]^{\E^d}.
\end{equation}
Here and below $\lambda\in(0,1)$ denotes the ellipticity ratio, which is fixed  throughout the paper. 
This allows, for instance, to interpret $\nabla^*a\nabla$ as either the operator of a ``conductance model'' (i.e.\ the solution of~\eqref{PDE} is a potential on a network of resistors) or the generator of a random walk on $\Z^d$ with jump rates across edges described by $a$.
Note that if we interpreted $\nabla^* a \nabla$ as a discretization of a continuum operator $-\nabla \cdot a \nabla$, the coefficient field $a\in\Omega$ would be diagonal next to being symmetric and uniformly elliptic. In the discrete setting, diagonality is
known to be a sufficient (but not necessary) condition for the maximum principle to 
hold for $\nabla^* a\nabla u$. The maximum principle is a crucial ingredient for the estimates \eqref{P1.1} and \eqref{P1.2} on the
quenched Green function, on which our results rely.

\medskip

Our main object is the non-constant coefficient, elliptic, discrete Green function
$G(a;x,x')$ defined through $\nabla^*a\nabla G(a;\cdot,x')=\delta( \cdot -x')$, 
where $\delta$ stands for the discrete version of the Dirac distribution, i.e.\ 
\[
 \delta(x) = \begin{cases}
              1 \quad&\text{for } x = 0\\
	      0 \quad&\text{otherwise}
             \end{cases}\Bigg\}.
\]

We usually drop the argument $a$ and just write $G(x,y)$. Often, it is more convenient to appeal to the distributional
characterization:
\begin{equation}\label{Green_dist}
\forall\;\zeta(x):\quad\sum_e\nabla\zeta(e) a(e)\nabla G(e,x')=\zeta(x').
\end{equation}

Here and throughout the paper, derivatives are understood to fall on the edge variable.
We will always work in dimension $d\ge 2$. Dimension $d=2$ needs a bit more care in terms of the
definition of the Green function. Since we are only interested in {\it gradient} estimates,
this is merely technical and will be ignored here. Sometimes, it is more convenient to think of $\nabla G$ as the limit of $\nabla G_T$ as $T\to\infty$ where $G_T$ is the Green's function with a massive term in the sense that
\begin{equation}\label{Green_massive}
 T^{-1} G_T(\cdot,x') + \nabla^* a \nabla G_T(\cdot,x') = \delta(\cdot-x');
\end{equation}
this is the case in the proof of Proposition \ref{P}. At other times,
it is more convenient to think in terms of an approximation via periodization in the sense of
\begin{equation}\label{Green_periodic}
 \nabla^* a(x) \nabla G_L(\cdot,x') = \sum_{z\in\Z^d} \delta(\cdot-x'-Lz) - L^{-d};
\end{equation}
this is the case in the proof of Lemma \ref{L3}.


\section{Assumptions on the ensemble and main result}\label{sec:main}

We are given a probability measure on the space $\Omega$ of uniformly elliptic, diagonal coefficient fields
(endowed with the product topology), cf.\ \eqref{Omega} in the previous section. Following the convention in
statistical mechanics, we call this probability measure an \emph{ensemble} and denote the associated 
ensemble average (i.e.\ the expected value) by $\expec{\cdot}$. Functions $\zeta:\Omega\to\R$ will also be called \emph{random variables}.
Note that $\Z^d$ acts on $\E^d$ by translation and we denote by $b+x\in\E^d$ the edge $b\in\E^d$ shifted by $x\in\Z^d$.
With this definition, we assume that $\expec{\cdot}$ is \emph{stationary} in the sense that for any shift vector $z\in\mathbb{Z}^d$,
the shifted coefficient field $a(\cdot+z):=(\E^d\ni e\mapsto a(e+z))\in\Omega$ has the same distribution
as $a$. We also note that the Green function is shift-invariant or
stationary in the sense that $G(a(\cdot+z);x,y)=G(a;x+z,y+z)$.

\medskip

Besides stationarity, the main assumption on the ensemble of coefficients and only probabilistic tool will be a variant of the 
logarithmic Sobolev inequality (LSI). It constitutes a quantification of ergodicity. In Section \ref{sec:LSI}, we will
comment on the LSI and the related spectral gap inequality --- there we will also describe the relation between this LSI and the usual LSI.

\begin{definition}\label{Dcont}{\rm [Logarithmic Sobolev inequality].}
Let $\langle\cdot\rangle$ be a (not necessarily stationary) ensemble of coefficients $a$.

We say $\langle\cdot\rangle$ satisfies a logarithmic Sobolev inequality (LSI) with constant $\rho>0$ if
for all random variables {$\zeta:\Omega \to \R$, we have that
\begin{equation}\label{LSI}
\Big\langle \zeta^2\log\frac{\zeta^2}{\langle \zeta^2\rangle}\Big\rangle
\le\frac{1}{2\rho} \Big\langle\sum_{e\in\E^d} \Big( \osc_{a(e)} \zeta \Big)^2\Big\rangle,
\end{equation}
where the oscillation is to be taken over all values of $a(e)\in[\lambda,1]$, i.e.\
over all coefficient fields $\tilde a \in \Omega$ that coincide with $a$ outside of $e\in\E^d$ (i.e.\ $\tilde a(b) = a(b)$ for all $b\neq e$). In formulas:
\begin{multline*}
 \Big( \osc_{a(e)} \zeta \Big)(a) = \sup\{ \zeta(\tilde a) \ | \ \tilde a\in \Omega \text{ s.t.\ } \tilde a(b)=a(b)\ \forall b\neq e \}\\ - \inf\{ 
\zeta(\tilde a) \ | \ \tilde a\in \Omega \text{ s.t.\ } \tilde a(b)=a(b)\ \forall b\neq e  \}.
\end{multline*}}
\end{definition}

Note that the difference between the LSI \eqref{LSI} and the usual LSI, see \eqref{LSI_usual}, lies in the use of the oscillation instead of the partial derivative $\frac{\partial \zeta}{\partial a(e)}$. The merit of this form is that it is satisfied by \emph{any} ensemble of independent, identically distributed coefficients $(a(e))_{e\in\E^d}$, cf.\ Lemma \ref{LSI_iid} below.
Our main result is: 

\medskip

\begin{theorem}\label{T}
Let $\langle\cdot\rangle$ be stationary and satisfy the LSI \eqref{LSI} with constant $\rho>0$, see Definition
\ref{Dcont}.
Then for all $1\le p<\infty$, $x\in\Z^d$ and $b,b'\in\mathbb{E}^d$, we have that
\begin{align}
\langle|\nabla\nabla G(b,b')|^{2p}\rangle^\frac{1}{2p} &\le C(d,\lambda,\rho,p)(|b-b'|+1)^{-d},\label{T.1}\\
\langle|\nabla G(b,x)|^{2p}\rangle^\frac{1}{2p}&\le C(d,\lambda,\rho,p)(|b-x|+1)^{1-d}.\label{T.3}
\end{align}
\end{theorem}
We furthermore let $|b|$ denote the Euclidean distance of the midpoint of the edge $b$ from the origin and $|b-b'|$ the distance between the midpoints of the two edges $b$ and $b'$. Recall that $b+x$ denotes the edge $b$ shifted by $x$. Here and in the sequel, $C(d,\lambda,\rho,p)$ stands for a generic constant that only depends on dimension $d\ge 2$, on the ellipticity ratio $\lambda>0$,
on the LSI constant $\rho>0$, and on the exponent of integrability $p<\infty$.

\smallskip

We defer the proof of Theorem \ref{T} until Subsection \ref{ssec:T}.
Clearly, the spatial decay rates in Theorem \ref{T} are optimal, since those are the decay
rates of the constant coefficient Green function, see for instance \cite[Theorem 4.3.1]{LawlerLimic}. 
Note that we may assume without loss of generality that $x=0$ in \eqref{T.3} since stationarity of $\langle\cdot\rangle$ and $G$ implies
\[\langle |\nabla G(a;b,x)|^{2p}\rangle
=\langle |\nabla G(a(\cdot-x);b,x)|^{2p}\rangle
=\langle|\nabla G(a;b-x,0)|^{2p} \rangle.\]
An interesting aspect of Theorem \ref{T} is the following: The \emph{quenched} versions of
(\ref{T.1}) and (\ref{T.3}) are false, i.e.\ the {\it uniform} in $a$ and \emph{point-wise} in $x$
estimates $|\nabla\nabla G(a;e,b)|\le C(d,\lambda)(|e-b|+1)^{-d}$ 
and $|\nabla G(a;e,0)|\le C(d,\lambda)(|e|+1)^{d-1}$ do \emph{not} hold (while suitably \emph{spatially} averaged
versions of both estimates do hold uniformly in $a$); see our discussion in Section \ref{sec:DeGiorgi} below.

\bigskip

An easy consequence is the following generalized variance estimate on $G$ itself:

\begin{corollary}\label{C}
Let $\langle\cdot\rangle$ be as in Theorem \ref{T}. Then we have that
\begin{equation}\label{C.1}
\Big\langle \big| G(x,0) - \langle G(x,0) \rangle \big|^{2p} \Big\rangle^{\frac{1}{p}} \le C(d,\lambda,\rho,p)
                                                                            \begin{cases}
                                                                             (|x|+1)^{2(1-d)} &d>2\\
									     (|x|+1)^{-2} \log(|x|+2)  &d=2
                                                                            \end{cases}\Bigg\}
\end{equation}
for all $x\in\mathbb{Z}^d$ and $1\le p < \infty$.
\end{corollary}

The proof of Corollary \ref{C} will be given in Subsection \ref{ssec:C}.

\begin{remark}
We note that the estimate in Corollary \ref{C} is optimal in the scaling of the spatial decay. This can
be seen by developing to leading order in a small ellipticity ratio $1-\lambda\ll 1$. We expand upon this argument (for
the special case of $p=1$) in Subsection \ref{ssec:opt_C} after the proof of Corollary \ref{C}.
\end{remark}


\section{Homogenization error}\label{sec:hom}

In the same vein as Corollary \ref{C}, Theorem \ref{T} allows to give optimal estimates on
the random part of the homogenization error. These extend the results by Conlon and Naddaf
\cite[Theorem 1.2, Theorem 1.3]{ConlonNaddaf} from small ellipticity ratio (i.e.\ $1-\lambda\ll 1$)
to arbitrary ellipticity ratio. For the ``strong error" (see below for an explanation of this wording) \cite[Theorem
1.2]{ConlonNaddaf} in $d>3$, this was already achieved by Gloria \cite[Theorem 2]{Gloria}. For all other cases, our
result appears to be new.
Let us be more precise: For a coefficient field $a:\E^d\to\R$ 
and a right-hand side $f:\Z^d\to\R$ we consider the solution $u:\Z^d\to \R$ of
\begin{equation}\label{u_eps}
 \nabla^*a\nabla u = f \quad\text{on $\Z^d$.}
\end{equation}
In order for (\ref{u_eps}) to have a unique solution that decays (i.\ e.\ $\lim_{|x|\to\infty}u(x)=0$), 
we assume for simplicity that $f$ is compactly supported (and furthermore is of zero spatial average in the case of
$d=2$).
By the random part of the homogenization error, we understand the ``fluctuations'' $u(x)-\langle u(x)\rangle$.
These are expected to be small (w.\ r.\ t.\ the size of $u(x)$ itself)
if $f(x)$ varies only slowly w.\ r.\ t.\ to the lattice spacing.
In our notation, the lattice spacing is unity, so that a natural model for a right-hand side that
has a large characteristic scale $L\gg 1$ is given by $f(x)=L^{-2}\hat f(\frac{x}{L})$ for some bounded
and compactly supported ``mask'' $\hat f(\hat x)$, $\hat x\in\mathbb{R}^d$. The scaling $L^{-2}$ of
the amplitude of $f$ is motivated as follows: In the rescaled variables $\hat x$, (\ref{u_eps}) now assumes the
suggestive form of
\begin{equation}\label{model}
 \big( \nabla_\epsilon^*a({\textstyle\frac{\cdot}{\epsilon}})\nabla_\epsilon u\big)({\hat x}) = \hat f({\hat x})  \quad \text{for all $\hat{x}\in\epsilon\Z^d$,}
\end{equation}
where $\epsilon:=L^{-1}$ is the ratio of the lattice spacing to the characteristic scale of the r.-h.\ s.
and where $\nabla_\epsilon$ denote finite differences for the rescaled lattice $\epsilon\mathbb{Z}^d$
(i.\ e.\ $\nabla_{\epsilon}u([\hat x, \hat x+\epsilon e_i])=\epsilon^{-1}(u(\hat x+\epsilon e_i)-u(\hat x))$).

\medskip

The size of the fluctuations will be measured in two different ways.
\begin{itemize}
\item Corollary \ref{Conlon_T1.2}: Here, the fluctuations will be controlled in a {\it strong} way in the sense that
we estimate the (discrete) $\ell^p(\mathbb{Z}^d)$-norm $\left(\sum_{x}|u-\langle
u\rangle|^p\right)^{1/p}$ 
of the fluctuations. This will be done for arbitrary stochastic moments (the role played by $rp$).
Corollary \ref{Conlon_T1.2} is the generalization of \cite[Theorem 1.2]{ConlonNaddaf} as well as \cite[Theorem 2]{Gloria}. For our model
right-hand side,
$f(x)=\epsilon^2\hat f(\epsilon x)$ with bounded and compactly supported $\hat f$,
the fluctuations are (up to a logarithmic correction for $d=2$) of the
order of $\epsilon$ in this measure, see (\ref{C2.1eps}).
\item Corollary \ref{Conlon_T1.3}: Here, the fluctuations will be controlled in a {\it weak} way in the sense
that we only estimate {\it spatial averages} $\sum_{x}(u-\langle u\rangle) g$ of the fluctuations,
with {\it deterministic} averaging function $g(x)$. Again, this will be done for arbitrary stochastic moments (the role
played by $r$).
Corollary \ref{Conlon_T1.3} is the generalization of \cite[Theorem 1.3]{ConlonNaddaf}. For our model right-hand side
$f(x)=\epsilon^2\hat f(\epsilon x)$ with bounded and compactly supported $\hat f$,
and an averaging function of the form $g(x)=\hat g(\epsilon x)$  with bounded and compactly supported $\hat g$,
the fluctuations are $O(\epsilon^{d/2})$ in this measure, see (\ref{C3.1eps}). (Here, there is no logarithmic correction even for $d=2$.)
\end{itemize}

\begin{corollary}\label{Conlon_T1.2}
Let $\langle\cdot\rangle$ be as in Theorem \ref{T}; 
for compactly supported right-hand side $f(x)$, consider the decaying solution $u(x)$ to \eqref{u_eps}. 
Let the spatial integrability exponents $2\le p<\infty$ and $1 < q < \infty$
be related through $\frac{1}{q}=\frac{1}{d}+\frac{1}{p}$. 

\smallskip

In case of $d>2$, we have for all $r<\infty$:
\begin{equation}\label{C2.1}
\bigg\langle \bigg(\sum_{x}\big|u-\langle u\rangle\big|^{p}
\bigg)^r \bigg\rangle^{\frac{1}{rp}} \le
C(d,\lambda,\rho,p,r) \Big(\sum_{x}|f|^q \Big)^{\frac{1}{q}}.
\end{equation}

\smallskip

In case of $d=2$, we additionally require $p>2$ (so that $q>1$) and that $f$ is supported in $\{x:|x|\le R\}$ for some $R\ge 1$. Then we have for all $r<\infty$:
\begin{equation}\label{C2.2}
\Bigg\langle\Bigg(\sum_{x:|x|\le R}\big|u-\langle u\rangle\big|^{p}\Bigg)^r\Bigg\rangle^{\frac{1}{rp}}
 \le C(\lambda,\rho,p,r)\,(\log^{\frac{1}{2}}R)\, \Big(\sum_{x}|f|^q \Big)^{\frac{1}{q}}.
\end{equation}

\end{corollary}

\begin{corollary}\label{Conlon_T1.3}
Let $\langle\cdot\rangle$ be as in Theorem \ref{T}; 
for compactly supported right-hand side $f(x)$, consider the decaying solution $u(x)$ to \eqref{u_eps}. 
Let the averaging function $g(x)$ be compactly supported. Let the two integrability exponents
{$1<q,\tilde q < d$} be related by $\frac{1}{q}+\frac{1}{\tilde q}=\frac{2}{d}+\frac{1}{2}$.
Then we have for all $r<\infty$:
\begin{equation}\label{C3.1}
\bigg\langle\bigg|\sum_{x} (u-\langle u\rangle) g\bigg|^{r} \bigg\rangle^{\frac{1}{r}}\\
 \le C(d,\lambda,\rho,r) 
\bigg(\sum_{x} |f|^q\bigg)^{\frac{1}{q}} 
\bigg(\sum_{x} |g|^{\tilde q}\bigg)^{\frac{1}{\tilde q}}.
\end{equation}
\end{corollary}

Corollaries \ref{Conlon_T1.2} and \ref{Conlon_T1.3} will be proved in Subsections \ref{ssec:Conlon}.
For the convenience of the reader, we express the results of both corollaries in terms of
the rescaled variable $\hat x=\epsilon x$, the model right-hand side $f(x)=\epsilon^2\hat f(\epsilon x)$
and the model averaging function $g(x)=\epsilon^d\hat g(\epsilon x)$; we also rewrite the
solution itself in terms of $u(x)=\hat u_\epsilon(\epsilon x)$. In this notation,
(\ref{C2.1}) (multiplied by $\epsilon^{d/p}$) turns into
\begin{multline}\label{C2.1eps}
\bigg\langle \bigg(\epsilon^d\sum_{\hat x\in\epsilon\mathbb{Z}^d}\big|\hat u_\epsilon-\langle \hat
u_\epsilon\rangle\big|^{p}
\bigg)^r \bigg\rangle^{\frac{1}{rp}}\\
\le
C(d,\lambda,\rho,p,r) \epsilon \Big(\epsilon^d\sum_{\hat x\in\epsilon\mathbb{Z}^d}|\hat f|^q \Big)^{\frac{1}{q}}
\;\le\;C(d,\lambda,\rho,r,\hat f) \epsilon.
\end{multline}
Note that this can be interpreted as the discrete version of
$$
\bigg\langle \bigg(\int_{\mathbb{R}^d}\big|\hat u_\epsilon-\langle \hat u_\epsilon\rangle\big|^{p}d\hat x
\bigg)^r \bigg\rangle^{\frac{1}{rp}}
\le
C(d,\lambda,\rho,p,r) \epsilon \Big(\int_{\mathbb{R}^d}|\hat f|^qd\hat x \Big)^{\frac{1}{q}},
$$
which highlights the $O(\epsilon)$-nature of the ``spatially strong'' error.

\smallskip

Likewise, (\ref{C3.1}) turns into
\begin{align}\label{C3.1eps}
&\bigg\langle\bigg|\epsilon^d\sum_{\hat x\in\epsilon\mathbb{Z}^d} 
(\hat u_\epsilon-\langle \hat u_\epsilon\rangle) \hat g\bigg|^{r} \bigg\rangle^{\frac{1}{r}}\nonumber\\
&\le C(d,\lambda,\rho,r) \epsilon^\frac{d}{2}
\bigg(\epsilon^d\sum_{\hat x\in\epsilon\mathbb{Z}^d} |\hat f|^q\bigg)^{\frac{1}{q}} 
\bigg(\epsilon^d\sum_{\hat x\in\epsilon\mathbb{Z}^d} |\hat g|^{\tilde q}\bigg)^{\frac{1}{\tilde q}}\nonumber\\
&\le C(d,\lambda,\rho,r,\hat f,\hat g) \epsilon^\frac{d}{2}.
\end{align}
As above, this can be seen as the discrete version of
\begin{equation}\nonumber
\bigg\langle\bigg|\int_{\mathbb{R}^d} 
(\hat u_\epsilon-\langle \hat u_\epsilon\rangle) \hat g\bigg|^{r}d\hat x \bigg\rangle^{\frac{1}{r}}
 \le C(d,\lambda,\rho,r) \epsilon^\frac{d}{2}
\bigg(\int_{\mathbb{R}^d} |\hat f|^qd\hat x\bigg)^{\frac{1}{q}} 
\bigg(\int_{\mathbb{R}^d} |\hat g|^{\tilde q}d\hat x\bigg)^{\frac{1}{\tilde q}},
\end{equation}
uncovering the $O(\epsilon^{d/2})$-nature of the ``spatially weak'' error.

\medskip

Let us make a couple of further more detailed remarks related to Corollaries \ref{Conlon_T1.2}
and \ref{Conlon_T1.3}. In case of Corollary \ref{Conlon_T1.2} and $d=2$,
we can use H\"older's inequality to establish an estimate also for $p=2$. However, in that case we
pay the price of an arbitrarily small power of $R$ on the right-hand side of (\ref{C2.2}).
We also note that the requirement that $f$ has compact support and that $u$ decays can be weakened: All we need 
is the Green function representation $u(x)=\sum_{y}G(x,y)f(y)$.
We conclude by pointing out that our argument does not require any
smoothness assumptions on $\hat f(\hat x)$ and $\hat g(\hat x)$ beyond (uniform) boundedness to obtain (\ref{C2.1eps}) and (\ref{C3.1eps}).

\medskip

The central limit theorem (CLT) scaling $O(\epsilon^{d/2})$ of the weak error 
seems to suggests that $u_\epsilon(x)$ behaves like a 
random field of amplitude $O(1)$ and integrable correlations. In fact, this is misleading,
as can be seen by distinguishing the scale $\frac{1}{\epsilon}$ on which $f$ varies
from the scale $1\ll\frac{1}{\delta}\ll\frac{1}{\epsilon}$ on which we take the spatial average
with help of the function $g$. If Corollary~\ref{Conlon_T1.3} were true in the limiting case of $q=d$
(which is not the case since the Hardy-Littlewood-Sobolev inequality in Step~3 in the proof of Corollary~\ref{Conlon_T1.3} requires $q<d$), we would obtain
\begin{equation}\nonumber
\bigg\langle\bigg|\delta^d\sum_{x\in\mathbb{Z}^d} 
(\hat u_\epsilon (\epsilon x) - \langle \hat u_\epsilon(\epsilon x) \rangle) \hat g(\delta x)\bigg|^{r} \bigg\rangle^{\frac{1}{r}}\nonumber\\
\le C(d,\lambda,\rho,r,\hat f,\hat g) \epsilon \delta^{\frac{d}{2}-1}.
\end{equation}
This refined estimate does suggest that $\hat u_\epsilon(\epsilon x)$ behaves like a
random field of amplitude $O(\epsilon)$ and correlations that decay like the Green's function:
\begin{equation}\nonumber
 \big|\big\langle (\hat u_\epsilon(\epsilon x) - \langle \hat u_\epsilon(\epsilon x) \rangle) (\hat u_\epsilon(\epsilon y) - \langle \hat u_\epsilon(\epsilon y) \rangle) \big\rangle \big| \le C(d,\lambda,\rho,\hat f) \epsilon^2 (|x-y|+1)^{2-d}
\end{equation}
for all $x,y\in\Z^d$.
This scaling is natural, since it would follow from the (higher-order, two-scale) expansion $\hat u_\epsilon(\hat x) \approx u_{\mathrm{hom}}(\hat x) + \epsilon\sum_{k=1}^d\phi_k(\frac{\hat x}{\epsilon})\partial_k u_{\mathrm{hom}}(\hat x)$  and the expected --- but unproven --- estimate on the covariance of this corrector:
\begin{equation}\nonumber
 |\langle \phi_k(x) \phi_k(y) \rangle| \le C(d,\lambda,\rho) (|x-y|+1)^{2-d}
\end{equation}
for all $x,y\in\Z^d$.
In the above, the function $\phi_k$ is the corrector in direction $e_k$ (which is an $a$-harmonic function of affine behavior on large scales) and $u_{\mathrm{hom}}$ is the solution to the elliptic equation with homogenized coefficients.
We remark here that the above-mentioned expansion for $u_\epsilon$ was recently quantified by Gloria, Neukamm and the second author \cite{GloriaNeukammOtto2} using Theorem~\ref{T}. Indeed, there it is shown that the error in an $H^1$-norm in space and $L^2$-norm in probability for $u_\epsilon - u_{\mathrm{hom}} - \epsilon\sum_{k=1}^d\phi_k(\frac{\cdot}{\epsilon})\partial_k u_{\mathrm{hom}}$ is still of order $\epsilon$, cf.\ \eqref{C2.1eps}. In order to obtain this result, the authors also treat the so-called systematic error, which is the difference between $\langle u_\epsilon \rangle$ and $u_{\mathrm{hom}}$.

\medskip

A more traditional CLT-scaling has been established for the {\it energy density}.
For $g=f$, the weak measure of
fluctuations turns into a measure of fluctuations of the energy:
\[
\sum_{x}(u-\langle u\rangle)g=\sum_{e}a (\nabla u)^2-\Big\langle\sum_{e} a(\nabla u)^2\Big\rangle.
\]
If we set $u=\phi_k$, then the (stationary) energy density defines the homogenized diffusion coefficient. In \cite[Theorem 2.1]{GloriaOtto1},
it is shown that in the case of independent, identically distributed (i.\ i.\ d.) coefficients, the energy density of
the corrector has CLT scaling in the sense that spatial averages behave
as if the energy density was independent from site to site; in \cite[Proposition 7]{GloriaNeukammOtto}, that result has been generalized to ensembles that only satisfy a spectral gap condition. 
The scaling result has been substantially sharpened for i.\ i.\ d.\ ensembles:
In this situation, the fluctuations of the energy of the corrector
become more and more Gaussian as the box over which
the spatial average is taken increases. The latter result has been obtained by three different techniques:
Nolen \cite{Nolen} gives a quantitative estimate based on a differential characterization of Gaussian distributions
(second-order Poincar\'e inequality) and relies on the corrector estimates from \cite[Theorem 2.1]{GloriaOtto1}.
Biskup, Salvi, and Wolff \cite{Biskup} obtain
a more qualitative result using a Martingale decomposition of the spatially averaged
energy density (their result assumes small ellipticity
contrast $1-\lambda\ll 1$, but presumably could be extended using the results of \cite{GloriaNeukammOtto}).
Rossignol \cite{Rossignol} in turn uses an orthogonal decomposition of the space of coefficients (Walsh decomposition).


\section{Relation to De Giorgi's approach to elliptic regularity}\label{sec:DeGiorgi}

While our result heavily relies on the celebrated regularity theory for
scalar elliptic operators, connected with the names of De Giorgi, Nash, and Moser,
it also gives a new perspective on these results. We will specify the input from regularity
theory, namely Nash's (upper) bounds on the parabolic Green function, in the next section. We now 
address what we see as a new perspective on these results, namely on De Giorgi's result
on H\"older continuity of $a$-harmonic functions.

\medskip

An elementary consequence of the mean value property is the following Liouville principle: Harmonic functions that
grow sub-linearly must be constant. This holds for the constant-coefficient Laplacian
both on $\mathbb{R}^d$ and on $\mathbb{Z}^d$, but is no longer true for {\it variable}
coefficients, even if they are uniformly elliptic. Indeed, a well-known example \cite[Corollary 16.1.5]{Iwaniec} 
shows that for any $\alpha>0$, there exists an explicit coefficient field $\alpha^2\le a(z)\le 1$ such
that $u(z)={\mathcal Re}(|z|^{\alpha-1}z)$ is $a$-harmonic in $z\ni\mathbb{C}\cong\mathbb{R}^2$. 
We believe that this example can be adapted to the lattice $\mathbb{Z}^2$
(provided the condition of diagonality is relaxed to the condition that the
discrete maximum principle is valid, a setting to which our results presumably can be
extended).
A celebrated result of De Giorgi \cite[Theorem 2]{DeGiorgi} states that this is the worst-case scenario: For any
dimension $d$ and any ellipticity ratio $\lambda$, there exists an exponent $\alpha_0(d,\lambda)>0$
with the following property:
For any field of coefficients $\lambda\le a(x)\le 1$ and any
$a$-harmonic function $u(x)$, a bound of the form $|u(x)|\le C |x|^{\alpha_0}$ for $|x|\gg 1$ implies
that $u$ is constant. This result
holds both in $\mathbb{R}^d$ and in $\mathbb{Z}^d$ \cite[Proposition 6.2]{Delmotte}.
In this sense, while it is no longer true that ``sub-linear implies
constant'', it remains true that ``\emph{very} sub-linear implies constant''.

\medskip

De Giorgi's result is in fact more quantitative and can be rephrased as an inner regularity result in terms of H\"older continuity
with H\"older exponent $\alpha_0$:
For any harmonic function $u(x)$ on $\{x:|x|\le R\}$, the $\alpha_0$-H\"older modulus of continuity at zero is estimated by the supremum:
\begin{equation}\nonumber
\sup_{x:|x|\le R} \frac{|u(x)-u(0)|}{|x|^{\alpha_0}}\le C(d,\lambda)R^{-\alpha_0}\sup_{x:|x|\le R}|u(x)|.
\end{equation} 
To contrast De Giorgi's result with our result below, let us rephrase it as follows:
\begin{equation}\label{quenched_Holder}
\forall\;\lambda\le a(x)\le 1,\;\; \forall\; R < \infty:\quad
\sup_{u}\frac{\sup_{x:|x|\le R}\frac{|u(x)-u(0)|}{|x|^{\alpha_0}}}
{\frac{1}{R^{\alpha_0}}\sup_{x:|x|\le R}|u(x)|}
\le C(d,\lambda),
\end{equation}
where the outer supremum is taken over all $u(x)$ that satisfy $\nabla^*a\nabla u=0$ in $\{x:|x|\le R\}$.

\medskip

In this context, we will show in Subsection \ref{ssec:C2} that Theorem \ref{T} has the following Corollary.
\begin{corollary}\label{C2}
 For all $0 < \alpha < 1$, $p < \infty$, and $R < \infty$, we have
\begin{equation}\label{annealed_Holder}
 \Bigg\langle \Bigg(\sup_{u} \frac{\sup_{x:|x|\le R} \frac{|u(x) - u(0)|}{|x|^\alpha}}{\frac{1}{R^\alpha}
\sup_{x:|x|\le R} |u(x)|} \Bigg)^p \Bigg\rangle \le C(d,\lambda,\alpha,p),
\end{equation}
where the outer supremum is taken over all $u(x)$ that satisfy $\nabla^* a \nabla u = 0$ in $\{x:|x| \le R\}$.
\end{corollary}

\medskip

{Loosely speaking, Corollary \ref{C2} implies that for ``most'' coefficient fields,
an $a$-harmonic function $u(x)$ is H\"older continuous with an exponent
\emph{arbitrarily close to one}. More precisely, the modulus of near-Lipschitz continuity
of $u(x)$ in some large ball is estimated by its supremum in the concentric ball of twice the radius with a ``quenched'' constant $C(a)$ with all moments bounded independently of the radius. Indeed, with the same proof the numerator in Corollary \ref{C2} can be chosen as the full H\"older-norm on $\{x:|x| \le \frac{R}{2}\}$. Furthermore it is straight-forward to extend the result to functions $\nabla^* a \nabla u = f$ if we include the $\ell^d$-norm of $f$ over $\{x:|x|\le R\}$ in the denominator.
The quantitative result of Corollary \ref{C2} has the Liouville principle as an easy corollary: For almost
every $a$, any sub-linear $a$-harmonic function $u$ must be constant. 
However, surprisingly for us, the (qualitative) Liouville principle} holds \emph{without any assumption on the
ensemble} $\langle\cdot\rangle$ besides stationarity! This is established in a very inspiring paper \cite[Theorem 3]{Benjaminietal}. The
main ingredients for the short and elegant argument are
\begin{itemize}
\item The ``annealed'' estimate $\langle\sum_{x}|x|^2 G(t,x,0)\rangle\le C t$
on the second moments of the parabolic Green function $G(a;t,x,y)\stackrel{\mbox{\scriptsize short}}{=}G(t,x,y)$ 
(cf.\ \cite[(SBD)]{Benjaminietal}, see Subsection \ref{sec:ingredients} below for the definition of $G$),
which in our uniformly elliptic context even holds in its stronger ``quenched'' version, that is,
$\sum_{x}|x|^2 G(t,x,0)\le C t$.
\item The annealed estimate $-\langle\sum_{x}G(t,x,0)\log G(t,x,0)\rangle\le C \log t$
on the spatial entropy of the parabolic Green function $G$ (cf.\ \cite[p.12]{Benjaminietal}),
which in our context is an immediate consequence of the second moments estimate.
This ingredient is shown to imply the following annealed continuity property of $G$:
\[
 \bigg\langle \sum_{y}G(1,0,y)\sum_{x}\frac{|G(t,0,x)-G(t-1,y,x)|^2}{G(t,0,x)+G(t-1,y,x)} \bigg\rangle \le \frac{C}{t}
\]
for some sequence $t\to\infty$.
\end{itemize}
%


\section{Logarithmic Sobolev inequality}\label{sec:LSI}

In the following, we give a more detailed description of our use of the logarithmic Sobolev inequality and prove that
any i.\ i.\ d.\ ensemble satisfies Definition~\ref{Dcont}. LSI substitutes the
spectral gap inequality (SG) in prior work on quantitative stochastic homogenization.
SG has been introduced into the field by Naddaf and Spencer
\cite[Theorem 1]{NaddafSpencerunpublished} (in form of the Brascamp-Lieb inequality)
and used most recently in \cite[Lemma 2.3]{GloriaOtto1} in an indirect way
and in \cite{GloriaNeukammOtto} explicitly. The LSI follows like SG from
the property that there is an integrable fall-off of correlations in the sense of
a uniform mixing condition \`a la Dobrushin-Shlosman, see for instance
\cite[Theorem 1.8 c)]{StroockZegarlinski} for a discrete setting.
Both SG and LSI quantify ergodicity of the ensemble, see for instance the discussion in \cite[Chapter~4]{GloriaNeukammOtto}. Recall that the usual LSI in this setting (with continuum derivative) would read
\begin{equation}\label{LSI_usual}
 \Big\langle \zeta^2\log\frac{\zeta^2}{\langle \zeta^2\rangle}\Big\rangle
\le\frac{1}{2\rho}
\Big\langle\sum_{e\in\E^d} \Big(\frac{\partial \zeta}{\partial
a(e)}\Big)^2\Big\rangle.
\end{equation}
In the LSI of Definition \ref{Dcont}, we have simply changed the derivative by an oscillation in order to capture ensembles whose marginal distribution contains atoms, as we shall explain now.

\medskip

Both SG and LSI are based on the notion of a {\it vertical derivative} (here, the oscillation) that defines a Dirichlet
form and thus a reversible dynamics, namely Glauber dynamics, on the space of coefficient fields
(the word ``vertical'' is used
to distinguish this derivative from the ``horizontal'' derivative naturally arising in
stochastic homogenization, but not used in this paper). In the earlier work
on stochastic homogenization and motivated by field theories, see \cite{NaddafSpencer},
the version of SG that is based on the {\it continuum}
vertical derivative (as on the r.h.s.\ of \eqref{LSI_usual}) has been used \cite{NaddafSpencerunpublished}. However, this assumption 
rules out the natural example of coefficients with a single-site distribution that only
assumes a {\it finite} number of values (Bernoulli). Hence
in order to treat arbitrary single-site distributions, we are
forced to consider the version of LSI found in Definition \ref{Dcont}. A SG inequality based on the
oscillation was already considered in \cite[Lemma 2.3]{GloriaOtto1}.

\medskip

The LSI has been of great use in the setting of stochastic processes and diffusion semi-groups, for the first time
introduced in generality by Gross \cite{Gross}. It implies SG as well as concentration of measure \cite[Chapter 5]{Ledoux} and is equivalent to the notion of
hyper-contractivity, see \cite[Theorem 1]{Gross} or \cite[Theorem 4.1]{GuionnetZegarlinski}. Incidentally, hyper-contractivity was
first observed in the Gaussian context by Nelson \cite{Nelson1}, see
\cite{Nelson2} for an improved result. It is thus the older notion and in fact motivated the (somewhat implicit)
introduction of LSI by Federbush~\cite{Federbush}. We refer to \cite{GuionnetZegarlinski} for a recent exposition on
LSI.

\medskip

The result of this section is that any independent, identically distributed coefficient-field satisfies the LSI
\eqref{LSI} of Definition \ref{Dcont}.
\begin{lemma}\label{LSI_iid}
 Consider an ensemble $\langle\cdot\rangle$ of i.\ i.\ d.\ coefficients on each edge with arbitrary marginal distribution on $[\lambda,1]$.
Then \eqref{LSI} holds, i.e.
\[
\Big\langle \zeta^2\log\frac{\zeta^2}{\langle \zeta^2\rangle}\Big\rangle
\le\frac{1}{2\rho}
\Big\langle\sum_{e\in\E^d}\Big(\osc_{a(e)\in[\lambda,1]}\zeta \Big)^2\Big\rangle
\]
for all functions $\zeta$ of the coefficient field $a$. The constant $\rho$ may be taken
to be $\rho = \frac{1}{8}$.
\end{lemma}
Lemma \ref{LSI_iid} is an immediate consequence of the following two lemmas. The first one shows that any single-edge distribution on
$[\lambda,1]$ satisfies the LSI in Definition \ref{Dcont}.
\begin{lemma}\label{L.LSI_ss}
 Let $\langle \cdot \rangle$ be any distribution on $[\lambda,1]$. Then we have that
\begin{equation}\label{LSI_ss}
 \Big\langle \zeta^2\log\frac{\zeta^2}{\langle \zeta^2\rangle}\Big\rangle
\le\frac{1}{2\rho} \Big(\osc_{a\in[\lambda,1]}\zeta \Big)^2
\end{equation}
for all functions $\zeta:[\lambda,1] \to \R$. In fact, the constant $\rho = \frac{1}{8}$ will do.
\end{lemma}
The next lemma shows that the LSI in Definition \ref{Dcont} satisfies the tensorization principle.
\begin{lemma}\label{L.LSI_tensor}
 Let $\langle \cdot \rangle$ be an ensemble consisting of independent distributions on the edges such that each
single-edge distribution satisfies the LSI \eqref{LSI_ss} with the same constant $\rho$. Then $\langle \cdot
\rangle$ itself satisfies the LSI \eqref{LSI} with constant $\rho$.
\end{lemma}
The proofs of Lemmas \ref{L.LSI_ss} and \ref{L.LSI_tensor} will be given in Subsection \ref{ssec:L.LSI}.


\section{Main ingredients of the proof}\label{sec:ingredients}

Loosely speaking, our approach consists in upgrading the 
(optimal) annealed estimates of Delmotte and Deuschel
\cite[Theorem 1.1]{DeuschelDelmotte} in terms of the integrability exponent $p$.

\begin{proposition}\label{P}{\rm [Delmotte and Deuschel]}.
Let $\langle\cdot\rangle$ be stationary. Then we have for all $b,b'\in\mathbb{E}^d$ and $x\in\Z^d$:
\begin{align}
\langle|\nabla\nabla G(b,b')|\rangle&\le C(d,\lambda)(|b-b'|+1)^{-d},\label{T.2}\\
\langle|\nabla G(b,x)|\rangle&\le C(d,\lambda)(|b-x|+1)^{1-d}.\label{I.1}
\end{align}
\end{proposition}

More precisely, we refer to the estimates (1.4) and (1.5a) in \cite[Theorem 1.1]{DeuschelDelmotte}
on the discrete {\it parabolic} Green function $G(t,x,y)=G(a;t,x,y)$ (i.e.\ the solution of
$\partial_tG(t,\cdot,y)+\nabla^*a\nabla G(t,\cdot,y)=0$ with $G(t=0,x,y)=\delta(x-y)$) that in our notation imply
for any weight exponent $\alpha<\infty$:
\begin{align}
\langle|\nabla\nabla G(t,b,b')|\rangle&\le C(d,\lambda,\alpha)(t+1)^{-\frac{d}{2}-1}
\Big(\frac{|b-b'|^2}{t+1}+1\Big)^{-\frac{\alpha}{2}},\label{P1.1}\\
\langle|\nabla G(t,b,x)|\rangle&\le C(d,\lambda,\alpha) (t+1)^{-\frac{d}{2}-\frac{1}{2}}
\Big(\frac{|b-x|^2}{t+1}+1\Big)^{-\frac{\alpha}{2}}.\label{P1.2}
\end{align}
(In fact, \cite{DeuschelDelmotte} establishes (\ref{P1.1}) and (\ref{P1.2})
with exponentially decaying weights instead of just algebraically decaying ones.)
Since the elliptic Green function can be inferred from the parabolic one
via $G(x,y)=\int_0^\infty G(t,x,y)dt$, these estimates imply (\ref{T.2}) and (\ref{I.1})
(by fixing some $\alpha>d$ and performing the change of variables $\hat t=|x|^{-2}(t+1)$).
Actually, \cite{DeuschelDelmotte} establishes (\ref{P1.2}) and thus (\ref{I.1}) in the stronger form 
where the $L^1$-norm $\langle|\cdot|\rangle$ is replaced by the $L^2$-norm 
$\langle|\cdot|^2\rangle^{1/2}$:
$\langle|\nabla G(t,b,x)|^2\rangle^{1/2}\le C(d,\lambda,\alpha)(t+1)^{-(d+1)/2}
(\frac{|b-x|^2}{t+1}+1)^{-\alpha/2}$.

\medskip

Let us point out
that the \emph{spatially point-wise annealed} estimates (\ref{P1.1}) and (\ref{P1.2}) are consequences of the
following \emph{spatially averaged quenched} estimates
\begin{align}
\sum_{x\in\Z^d}\bigg(\Big(\frac{|x|^2}{t+1}+1\Big)^\frac{\alpha}{2} G(t,x,0)\bigg)^2
&\le C(d,\lambda,\alpha)(t+1)^{-\frac{d}{2}},\label{P.4}\\
\sum_{b\in\E^d}\bigg(\Big(\frac{|b|^2}{t+1}+1\Big)^\frac{\alpha}{2}|\nabla G(t,b,0)|\bigg)^2
&\le C(d,\lambda,\alpha)(t+1)^{-\frac{d}{2}-1}\label{P.3}.
\end{align}
The first estimate \eqref{P.4} is the (upper, off-diagonal part of the) celebrated Nash estimate \cite[Appendix]{Nash}. 
The discrete case was treated in full generality in \cite[Corollary 3.28]{CarlenKusuokaStroock}. 
The second estimate \eqref{P.3} is a consequence of the first one. For an elementary proof of both, 
we refer to \cite[Lemmas 24 and 25]{GloriaNeukammOtto}, with the
Nash inequality as only noteworthy ingredient. 
Let us point out how (\ref{P.3}) is implies (\ref{P1.1}). Using the semi-group property in form of $\nabla\nabla G(t,b,b')=\sum_{y}\nabla
G(\frac{t}{2},b,y) \nabla G(\frac{t}{2},y,b')$ we obtain by the triangle inequality for the weight, Cauchy Schwarz
in $\sum_{y}$ and the symmetry of $G(t,x,y)$ in $x$ and $y$:
\begin{align*}
 &\Big(\frac{|b-b'|^2}{t+1}+1\Big)^\frac{\alpha}{2}|\nabla\nabla G(t,b,b')|\\
 &\le \sum_{y\in\Z^d}\Big(\frac{2|b-y|^2}{t+1}+1\Big)^\frac{\alpha}{2}|\nabla G({\textstyle\frac{t}{2}},b,y)|
\Big(\frac{2|b'-y|^2}{t+1}+1\Big)^\frac{\alpha}{2}|\nabla G({\textstyle\frac{t}{2}},y,b')|\\
 &\le \Bigg( \sum_{y\in\Z^d}\bigg(\Big(\frac{2|b-y|^2}{t+1}+1\Big)^\frac{\alpha}{2}|\nabla G({\textstyle\frac{t}{2}},b,y)|\bigg)^2\\
 &\qquad\times\sum_{y\in\Z^d}\bigg(\Big(\frac{2|b'-y|^2}{t+1}+1\Big)^\frac{\alpha}{2}|\nabla
G({\textstyle\frac{t}{2}},b',y)|\bigg)^2\Bigg)^\frac{1}{2}.
\end{align*}
Note that the right-hand side of the last inequality does not allow for application of \eqref{P.3}, since the sum
is not in the variable in which the derivative is taken.
However, we take the expectation,
use the Cauchy-Schwarz inequality in $\langle\cdot\rangle$ and stationarity and symmetry in form of $\langle|\nabla
G({\textstyle\frac{t}{2}},b,y)|^2\rangle = \langle|\nabla G({\textstyle\frac{t}{2}},b-y,0)|^2\rangle$ to obtain
\begin{align*}
 &\Big(\frac{|b-b'|^2}{t+1}+1\Big)^\frac{\alpha}{2}\langle|\nabla\nabla G(t,b,b')|\rangle\\
 &\le \bigg( \sum_{y\in\Z^d}\Big(\frac{2|b-y|^2}{t+1}+1\Big)^{\alpha}\langle|\nabla G({\textstyle\frac{t}{2}},b,y)|^2\rangle\\
 &\qquad\times\sum_{y\in\Z^d}\Big(\frac{2|b'-y|^2}{t+1}+1\Big)^{\alpha}\langle|\nabla
G({\textstyle\frac{t}{2}},b',y)|^2\rangle\bigg)^\frac{1}{2}\\
 &\le\bigg( \Big\langle\sum_{y\in\Z^d}\Big(\frac{2|b-y|^2}{t+1}+2\Big)^{\alpha}|\nabla
G({\textstyle\frac{t}{2}},b-y,0)|^2\Big\rangle\\
 &\qquad\times\Big\langle\sum_{y\in\Z^d}\Big(\frac{2|b'-y|^2}{t+1}+2\Big)^{\alpha}|\nabla
G({\textstyle\frac{t}{2}},b'-y,0)|^2\Big\rangle
\bigg)^\frac{1}{2}.
\end{align*}
We now see that (\ref{P.3}) implies (\ref{P1.1}). The estimate (\ref{P1.2}) is derived via
the semi-group property in form of $\nabla G(t,b,x)=\sum_{y}\nabla G(\frac{t}{2},b,y)
G(\frac{t}{2},y,x)$ from the combination of (\ref{P.4}) and (\ref{P.3}) by an analogous argument.

\bigskip

Note that the estimates of Proposition~\ref{P} make \emph{no assumptions
on the ensemble besides stationarity}.
In order to pass from Proposition~\ref{P} to Theorem~\ref{T},
we need the assumption on the ensemble from Definition~\ref{Dcont}.
In fact, LSI enters only through the following lemma, which we shall prove in Subsection~\ref{ssec:L1}.
\begin{lemma}\label{L1}
Let $\langle\cdot\rangle$ satisfy LSI \eqref{LSI} with constant $\rho>0$. Then for
arbitrary $\delta>0$ and $1\le p < \infty$ and for any $\zeta:\Omega\to\R$, we have that
\begin{equation}
\langle|\zeta|^{2p}\rangle^\frac{1}{2p}
\le C(d,\rho,p,\delta)\langle|\zeta|\rangle
+\delta\Big\langle\bigg(\sum_e\Big(\osc_{a(e)} \zeta \Big)^2\bigg)^p\Big\rangle^\frac{1}{2p}.
\label{L1.9}
\end{equation}
\end{lemma}

The preceding lemma may be seen as a reverse H\"older inequality in probability: If one controls a bit (recall that $\delta > 0$ may be arbitrarily small) of the vertical derivative of a random variable $\zeta$, then its $L^1_{\langle \cdot \rangle}$-norm bounds its $L^{2p}_{\langle \cdot \rangle}$-norm. It can be seen as a softening of the concentration of measure phenomenon, which requires Lipschitz continuity of $\zeta$, cf.~\cite[Theorem~5.3]{Ledoux}.

\bigskip

In order to make use of Lemma \ref{L1}, we need to estimate the vertical derivatives of 
$\nabla\nabla G$ and $\nabla G$. The following lemma is at the core of our result.

\begin{lemma}\label{L2}
There exists an integrability exponent
$p_0=p_0(d,\lambda)<\infty$ such that for all $p\ge p_0$, {we have that
\begin{align}\nonumber
 &\sup_{b,b'\in\E^d}\Bigg\{(|b-b'|+1)^d\Big\langle\bigg(\sum_{e\in\E^d}\Big(\osc_{a(e)} \nabla\nabla G(b,b')
\Big)^{2}\bigg)^p\Big\rangle^\frac{1}{2p}\Bigg\}\\
 &\le C(d,\lambda,p) \sup_{b,b'\in\E^d}\Big\{ (|b-b'|+1)^d \langle|\nabla\nabla G(b,b')|^{2p}\rangle^\frac{1}{2p} \Big\}\label{L2.25}
\end{align}
and
\begin{align}\nonumber
 &\sup_{b\in\E^d,x\in\Z^d}\Bigg\{(|b-x|+1)^{d-1}\Big\langle\bigg(\sum_{e\in\E^d}\Big(\osc_{a(e)} \nabla G(b,x)
\Big)^{2}\bigg)^p\Big\rangle^\frac{1}{2p}\Bigg\}\\
&\le C(d,\lambda,p)\bigg( \sup_{b\in\E^d,x\in\Z^d} \Big\{ (|b-x|+1)^{d-1} \langle|\nabla
G(b,x)|^{2p}\rangle^\frac{1}{2p}\Big\}\label{L2.24}\\
&\qquad\qquad\quad + \sup_{b,b'\in\E^d}\Big\{(|b-b'|+1)^{d}\langle|\nabla\nabla
G(b,b')|^{2p}\rangle^\frac{1}{2p}\Big\}\bigg).\nonumber
\end{align}
}
\end{lemma}

For the proof, we refer to Subsection \ref{ssec:L2}. Note that in contrast to Proposition~\ref{P}, here the only assumption on the ensemble is LSI \eqref{LSI} --- in particular, \emph{Lemmas
\ref{L1} and \ref{L2} do not require stationarity} and stationarity enters the proof of Theorem~\ref{T} only through Proposition~\ref{P}.
The formulation of Lemma \ref{L2} shows that with our method, we first have to
estimate the mixed {\it second}
derivatives $\langle|\nabla\nabla G(b,b')|^{2p}\rangle$ before we can tackle
the {\it first} derivatives $\langle|\nabla G(b,0)|^{2p}\rangle$. It also reveals
that it is necessary to estimate {\it high moments} $p\ge p_0$ in $\langle\cdot\rangle$
in order to estimate {\it moderately low moments} 
like the fourth moment $\langle|\nabla G(b,0)|^{4}\rangle$
that is needed in the proof of Corollary~\ref{C}.

\bigskip

The preceding lemma relies on the following suboptimal, but {\it quenched} estimates
on the (elliptic) Green function:

\begin{lemma}\label{L3}{\rm [Gloria and Otto]}
There exists an exponent $\alpha_0=\alpha_0(d,\lambda)>0$
such that for all $R>0$ and $b\in\E^d$, we have that
\begin{align}
 R^{2\alpha_0}\sum_{e:R\le |e-b|< 2R}|\nabla\nabla G(e,b)|^2&\le C(d,\lambda),\label{L2.17}\\
 \sum_{e:R\le |e|< 2R}|\nabla G(e,0)|^2&\le C(d,\lambda).\label{I.2}
\end{align}
\end{lemma}

The estimate (\ref{I.2}) was established in the stronger (dimensionally optimal) form of
$\sum_{R\le |e|< 2R}|\nabla G(e,0)|^2\lesssim R^{2-d}$ in \cite[Lemma 2.9]{GloriaOtto1};
in its weaker form of (\ref{I.2}), it is straight forward for $d>2$. The proof of
estimate (\ref{I.2}) in \cite{GloriaOtto1} in case of $d=2$ is subtle and relied on an adaptation of
\cite{Dolzmannetal}.
{In Subsection \ref{ssec:L3}}, we will give an elementary argument for the estimate (\ref{L2.17}) which 
we could not find in the literature. We remark that the proof presented here does not make use of the maximum
principle (directly or indirectly) and therefore is also applicable to the case of systems, which we intend to use in future work.

\begin{remark}
 We mention that with the same proof, one obtains a periodic version of Theorem \ref{T} (with constants uniform in $L$)
for the Green function defined in \eqref{Green_periodic}. In that case, one just replaces the Euclidean distance $|x|$
on $\Z^d$ by its periodic version $\mathrm{dist}(x,L\Z^d)$ on the torus $\R/L\Z^d$. The periodic version of Proposition
\ref{P} follows as above from the quenched spatially averaged estimates of \cite[Theorem~3(b)]{GloriaNeukammOtto}. The same is true in the
presence of a massive term, cf.~\eqref{Green_massive}.
\end{remark}


\section{Proofs}\label{sec:proofs}

\subsection{Proof of Lemma \ref{L1}}\label{ssec:L1}

{\bf Step 1}. Result for $p=1$. We claim that for any $\delta > 0$ and all $\zeta(a)$:
\begin{equation}\label{L1.7}
 \langle \zeta^2 \rangle^{\frac{1}{2}} \le \bigg( \exp\Big(\frac{2}{\rho\delta^2}\Big) + \frac{\rho 
\delta^2}{2e}\bigg) \langle |\zeta| \rangle + \delta
\Big\langle \sum_e \Big(\osc_{a(e)} \zeta \Big)^2 \Big\rangle^{\frac{1}{2}},
\end{equation}
where $\rho$ denote the constant in the LSI, see Definition \ref{Dcont}.
By homogeneity, we may assume $\langle\zeta^2\rangle = 1$.
For all real-valued $\zeta$ we have that
\[
 \zeta^2 \le \begin{cases}
            \exp(\frac{2}{\rho\delta^2}) |\zeta|& \quad\text{if $|\zeta| \le \exp{\frac{2}{\rho\delta^2}}$}\\
	    {\textstyle\frac{\rho\delta^2}{4}} \zeta^2 \log \zeta^2 &\quad\text{if
$|\zeta|\ge\exp{\frac{2}{\rho\delta^2}}$}
           \end{cases}\Bigg\}.
\]
Since $x\log x$ is bounded from below by $\frac{1}{e}$, we have that $\frac{2}{e} |\zeta| + \zeta^2 \log \zeta^2 \ge 0$
for all $\zeta$. It follows that
\[
 \zeta^2 \le \bigg( \exp\Big(\frac{2}{\rho\delta^2}\Big) + \frac{\rho \delta^2}{2e} \bigg) |\zeta| +
\frac{\rho\delta^2}{4} \zeta^2 \log \zeta^2.
\]
Hence taking the expectation $\langle \cdot \rangle$ yields
\[
 \langle \zeta^2 \rangle \le \bigg( \exp\Big(\frac{2}{\rho\delta^2}\Big) + \frac{\rho 
\delta^2}{2e} \bigg) \langle |\zeta| \rangle +
\frac{\rho\delta^2}{4}
\Big\langle \zeta^2 \log \zeta^2 \Big\rangle.
\]
Since $\langle\zeta^2\rangle = 1$, Young's inequality yields
\begin{align*}
 \langle |\zeta| \rangle &\le \frac{1}{2} \bigg( \exp\Big(\frac{2}{\rho\delta^2}\Big) + \frac{\rho 
\delta^2}{2e} \bigg) \langle
|\zeta| \rangle^2 + \frac{1}{2} \bigg( \exp\Big(\frac{2}{\rho\delta^2}\Big) + \frac{\rho 
\delta^2}{2e} \bigg)^{-1}\\
 &= \frac{1}{2} \bigg( \exp\Big(\frac{2}{\rho\delta^2}\Big) + \frac{\rho 
\delta^2}{2e} \bigg) \langle |\zeta| \rangle^2 + \frac{1}{2}
\bigg( \exp\Big(\frac{2}{\rho\delta^2}\Big) + \frac{\rho 
\delta^2}{2e} \bigg)^{-1} \langle \zeta^2 \rangle.
\end{align*}
Combining the last two estimates, we deduce
\[
 \langle \zeta^2 \rangle \le \bigg( \exp\Big(\frac{2}{\rho\delta^2}\Big) + \frac{\rho 
\delta^2}{2e} \bigg)^2 \langle |\zeta| \rangle^2
+ \frac{\rho\delta^2}{2} \Big\langle \zeta^2 \log \frac{\zeta^2}{\langle\zeta^2\rangle} \Big\rangle.
\]
Hence LSI yields
\[
 \langle \zeta^2 \rangle \le \bigg( \exp\Big(\frac{2}{\rho\delta^2}\Big) + \frac{\rho 
\delta^2}{2e} \bigg)^2 \langle |\zeta| \rangle^2
+ \delta^2 \Big\langle \sum_e \Big(\osc_{a(e)} \zeta \Big)^2 \Big\rangle
\]
and estimate \eqref{L1.7} follows from taking the square root and applying the inequality $\sqrt{\zeta+\xi} \le \sqrt{\zeta} + \sqrt{\xi}$ for all numbers $\zeta, \xi \ge 0$.

\medskip

{\bf Step 2}. We finish the proof of \eqref{L1.9}, i.e.\ show that
\[
\langle \zeta^{2p}\rangle^\frac{1}{2p}\le C(\rho,p,\delta)\langle|\zeta|\rangle
+\delta\Bigg(\Big\langle\bigg(\sum_e \Big(\osc_{a(e)} \zeta \Big)^2\bigg)^p\Big\rangle\Bigg)^\frac{1}{2p}
\]
for general $p \ge 1$.
To that end, we apply (\ref{L1.7}) to $\zeta$ replaced by $|\zeta|^p$: 
\begin{equation}\nonumber
\langle |\zeta|^{2p}\rangle\le C(\rho,p,\delta)\langle |\zeta|^{p}\rangle^2
+\delta\Big\langle \sum_e \Big(\osc_{a(e)} |\zeta|^p \Big)^2  \Big\rangle,
\end{equation}
where $C(\rho,p,\delta)$ denotes a generic constant only depending on $\rho$, $p$, and $\delta$.
Since $p<2p$, an application of H\"older's inequality in $\langle\cdot\rangle$
and Young's inequality on the first r.-h.\ s.\ term yields
\begin{equation}\label{L1.8}
\langle |\zeta|^{2p}\rangle\le C(\rho,p,\delta)\langle|\zeta|\rangle^{2p}
+2\delta\Big\langle \sum_e \Big(\osc_{a(e)} |\zeta|^p \Big)^2  \Big\rangle.
\end{equation}
Now we use that
\[
 \osc_{a(e)} |\zeta|^p \le C(p) \bigg( |\zeta|^{p-1} \osc_{a(e)} \zeta + \Big( \osc_{a(e)} \zeta \Big)^p \bigg)
\]
which follows from the elementary inequality $|\zeta^p - \xi^p| \le C(p) (\zeta^{p-1} |\zeta-\xi| + |\zeta-\xi|^p)$ for all numbers $\zeta,\xi > 0$ and the triangle inequality in form of $\osc_{a(e)} |\zeta| \le \osc_{a(e)} \zeta$. Hence \eqref{L1.8} yields
\begin{multline}\label{L1.osc^p}
\langle |\zeta|^{2p}\rangle\le C(\rho,p,\delta)\langle|\zeta|\rangle^{2p} + 2 C(p) \delta\Big\langle |\zeta|^{2p-2} \sum_e \Big(\osc_{a(e)}
\zeta \Big)^{2} \Big\rangle\\
+2 C(p) \delta\Big\langle \sum_e \Big(\osc_{a(e)} \zeta \Big)^{2p}  \Big\rangle.
\end{multline}
The last term on the right-hand side may be estimated by discreteness:
\begin{equation}\label{L1.discrete}
 \Big\langle \sum_e \Big(\osc_{a(e)} \zeta \Big)^{2p}  \Big\rangle \le \Big\langle \bigg( \sum_e \Big(\osc_{a(e)} \zeta \Big)^2
\bigg)^p \Big\rangle
\end{equation}
Furthermore, H\"older's inequality followed by Young's inequality yields
\begin{align}\nonumber
 \Big\langle |\zeta|^{2p-2} \sum_e \Big(\osc_{a(e)} \zeta \Big)^2 \Big\rangle &\le \langle |\zeta|^{2p} \rangle^{1-\frac{1}{p}} \Big\langle
\bigg( \sum_e \Big(\osc_{a(e)} \zeta \Big)^2 \bigg)^p \Big\rangle\\
 &\le \frac{1}{4 C(p) \delta} \langle |\zeta|^{2p} \rangle + (4C(p)\delta)^{p-1} \Big\langle \bigg( \sum_e \Big(\osc_{a(e)} \zeta \Big)^2
\bigg)^p \Big\rangle.\label{L1.Young}
\end{align}
Hence collecting \eqref{L1.osc^p}, \eqref{L1.discrete} and \eqref{L1.Young} yields
\begin{multline*}
 \langle |\zeta|^{2p}\rangle\le C(\rho,p,\delta)\langle|\zeta|\rangle^{2p} + 2\big(2 C(p) \delta + (4 C(p) \delta)^p \big) \Big\langle
\bigg( \sum_e \Big(\osc_{a(e)} \zeta \Big)^2
\bigg)^p \Big\rangle,
\end{multline*}
where we have absorbed the second term of \eqref{L1.Young} on the left-hand side.
By redefining $\delta$, we obtain (\ref{L1.9}).
%



\subsection{Proof of Lemma \ref{L3}}\label{ssec:L3}

We just give the proof of (\ref{L2.17}); for (\ref{I.2}), we refer to \cite[Lemma~2.9]{GloriaOtto1}. Note that in the stronger form $\sum_{e:R\le |b-e|< 2R}|\nabla\nabla G(e,b)|^2\le C(d,\lambda)R^{2-d-2\alpha_0}$, 
Estimate \eqref{L2.17} can also be seen
as a consequence of the following classical ingredients (which however would not hold in the systems case):
\begin{itemize}
\item the optimal decay of $G(x,y)$ itself, that is just needed in a spatially averaged sense of
$R^{-d}\sum_{y:R\le |x-y|< 2R} |G(x,y)-\bar G| \le C(d,\lambda)R^{2-d}$ (thanks to subtracting
the average $\bar G$ over the annulus $\{y:R\le|x-y|\le 2R\}$, this estimate also holds in $d=2$),
\item De Giorgi's H\"older continuity estimate, that then yields for some $\alpha_0=\alpha_0(d,\lambda)>0$ that
$\sup_{x:R\le |b-x|< 2R} |\nabla G(x,b)| \le C(d,\lambda)R^{2-d-\alpha_0}$,
\item Caccioppoli's estimate, that then yields
$\sum_{e:R\le |b-e|< 2R} |\nabla\nabla G(e,b)|^2 \le C(d,\lambda)
R^{2-d-2\alpha_0}$.
\end{itemize}

{\bf Step 1}. In this step, we derive the a priori estimate
\begin{equation}\label{L3.10}
\sum_{e}|\nabla\nabla G(e,b)|^2\le C(d,\lambda).
\end{equation}
Indeed, recall the weak formulation~\eqref{Green_dist} of the defining equation for $G$, i.e.\
\[
\forall\;\zeta(x):\quad\sum_{e}\nabla\zeta(e) a(e)\nabla G(e,x)=\zeta(x).
\]
Taking the derivative w.\ r.\ t.\ the variable $x$ along some edge $b$ yields
\begin{equation}\label{L3.18}
\forall\;\zeta(x):\quad\sum_{e}\nabla\zeta(e) a(e)\nabla\nabla G(e,b)=\nabla\zeta(b).
\end{equation}
The choice of $\zeta(x)=\nabla G(x,b)$ (we address the question of admissibility
of this test function below) yields
\[
\sum_{e} a(e)(\nabla\nabla G(e,b))^2
=\nabla\nabla G(b,b).
\]
Since $a(b)\ge\lambda$, this implies (\ref{L3.10}) in the explicit form of
\begin{equation}\label{L3.19}
\sum_{e}|\nabla\nabla G(e,b)|^2\le\lambda^{-2}.
\end{equation}

\medskip

We now turn to the question of admissibility of $\zeta(x)=\nabla G(x,b)$
as a test function for (\ref{L3.18}),
i.e.\ the question of decay as $|x|\to\infty$ of this function and its gradient.
This issue can be circumvented as in Step 3 below through approximation by the \emph{periodic} problem. 
More precisely, we consider the \emph{periodic} discrete elliptic Green function 
$G_L(x,x')=G_L(a,x,x')$ of period $L$. Up to additive constants, it is characterized by the weak equation
\begin{equation}\nonumber
\sum_{e\in\E^d\cap[-\frac{L}{2},\frac{L}{2})^d}\nabla\zeta(e) a(e)\nabla G_L(e,x')
=\zeta(x')-L^{-d}\sum_{x\in\Z^d\cap[-\frac{L}{2},\frac{L}{2})^d}\zeta(x)
\end{equation}
for all periodic $\zeta(x)$. Here $b\in\E^d\cap[-\frac{L}{2},\frac{L}{2})^d$ stands short for the set of all edges $b$ whose
midpoint is contained in the box $[-\frac{L}{2},\frac{L}{2})^d$. With the same argument as above, we obtain
\begin{equation}\label{energy_G_L}
\sum_{e\in\E^d\in[-\frac{L}{2},\frac{L}{2})^d}|\nabla\nabla G_L(e,b)|^2\le\lambda^{-2}.
\end{equation}
Since $G_L(x,x')$ converges point-wise to $G(x,x')$,
the latter implies (\ref{L3.19}) in the limit $L\to\infty$ by Fatou's lemma.
Incidentally, $\lim_{L\to\infty}\nabla G_L(e,x')$ may be taken as a 
definition of $\nabla G(e,x')$ in the case of $d=2$, 
where $G$ itself is not unambiguously defined.

\medskip


In the following steps, we use the fact that $u(x)=\nabla G(x,b)$ is $a$-harmonic away from the endpoints of $b$ to show that there
exists a decay exponent $\alpha_0(d,\lambda)>0$ such that for all $R\ge C(d)$ we have
\begin{equation}\label{L3.14}
\sum_{e:|b-e|\ge R}|\nabla u(e)|^2\le C(d,\lambda)R^{-2\alpha_0}\sum_{e:|b-e|\ge 1}|\nabla u(e)|^2.
\end{equation}
Together with (\ref{L3.10}), this implies (\ref{L2.17}). In Step 2, we will formally treat the
continuum whole-space case (where $a(x)$ is a uniformly elliptic matrix). In Step 3, we will show how to make the continuum case rigorous
by approximation through the continuum periodic case.
More precisely, using \eqref{L3.10}, we will directly prove the estimate \eqref{L2.17} in form of
\begin{equation}\label{L3.14b}
\sum_{e:|b-e|\ge R}|\nabla u(e)|^2\le C(d,\lambda)R^{-2\alpha_0}.
\end{equation}
In Step 4, we indicate the changes necessary to treat the discrete case.

\medskip

{\bf Step 2}. Formal derivation of the continuum version of (\ref{L3.14}), that is
\begin{equation}\label{L3.15}
\int_{\{x:|x|\ge R\}}|\nabla u|^2dx\le C(d,\lambda)R^{-2\alpha_0}\int_{\{x:|x|\ge 1\}}|\nabla u|^2 dx
\end{equation}
for $R\ge 1$ and a function $u(x)$ satisfying
\begin{equation}\label{L3.13}
-\nabla\cdot a\nabla u(x)=0\quad\mbox{in}\;\{x:|x|>1\}.
\end{equation}
Indeed, let $\eta(x)$ be a cut-off function for $\{x:|x|\ge 2R\}$ in $\{x:|x|\ge R\}$.
We test (\ref{L3.13}) with $\zeta=\eta^2 (u-\bar u)$, where $\bar u$ is the spatial
average of $u$ on the annulus $\{x:R\le |x|\le 2R\}$. It is a priori not clear that this
is an admissible test function for (\ref{L3.13}); we shall address this in the next step.
We appeal to the identity
\begin{equation}\label{L3.20}
\nabla(\eta^2 (u-\bar u))\cdot a\nabla u=\nabla(\eta (u-\bar u))\cdot a\nabla(\eta (u-\bar u))
-(u-\bar u)^2\nabla\eta\cdot a\nabla\eta,
\end{equation}
which in view of ellipticity in form of $\lambda |\xi|^2 \le \xi\cdot a(x) \xi \le |\xi|^2$ for all $\xi\in\R^d$ turns into the inequality
\begin{equation}\label{L3.21}
\nabla(\eta^2 (u-\bar u))\cdot a\nabla u\ge\lambda|\nabla(\eta (u-\bar u))|^2
-(u-\bar u)^2|\nabla\eta|^2.
\end{equation}
Hence from testing (\ref{L3.13}) we obtain 
\begin{equation}\nonumber
\lambda\int_{\R^d} |\nabla(\eta (u-\bar u))|^2dx\le\int_{\R^d} (u-\bar u)^2|\nabla\eta|^2dx,
\end{equation}
which by the choice of $\eta$ yields the Caccioppoli estimate
\begin{equation}\label{L3.22}
\int_{\{x:|x|\ge 2R\}}|\nabla u|^2dx\le C(d,\lambda)R^{-2}\int_{\{x:R\le |x|\le 2R\}}(u-\bar u)^2dx.
\end{equation}
By Poincar\'e's estimate on $\{x:R\le |x|\le 2R\}$ with mean value zero, this turns into
\begin{equation}\nonumber
\int_{\{x:|x|\ge 2R\}}|\nabla u|^2dx\le C(d,\lambda)\int_{\{x:R\le |x|\le 2R\}}|\nabla u|^2dx,
\end{equation}
which can be reformulated as
\begin{equation}\label{L3.2}
\int_{\{x:|x|\ge R\}}|\nabla u|^2dx\le C(d,\lambda)\int_{\{x:R\le |x|\le 2R\}}|\nabla u|^2dx.
\end{equation}

\medskip

A standard iteration argument now leads from \eqref{L3.2} to \eqref{L3.15}:
Introducing the notation $I_k:=\int_{\{x:|x|\ge 2^k\}}|\nabla u|^2dx$, estimate \eqref{L3.2} reads
\begin{equation}\nonumber
\forall k\in\{0,1,\cdots\}\quad I_{k}\le C(d,\lambda)(I_k-I_{k+1}),
\end{equation}
which with help of $\theta=\theta(d,\lambda):=1-\frac{1}{C}<1$ can be reformulated 
\begin{equation}\nonumber
\forall k\in\{0,1,\cdots\}\quad I_{k+1}\le \theta I_k,
\end{equation}
or with help of $\alpha_0=\alpha_0(d,\lambda):=\frac{-\log\theta}{2\log 2}$ as
\begin{equation}\nonumber
\forall k\in\{0,1,\cdots\}\quad I_{k}\le \theta^k I_0=(2^k)^{-2\alpha_0}I_0.
\end{equation}
In the original notation, this implies (\ref{L3.15}) in form of
\begin{equation}\nonumber
\forall R\ge 1\quad \int_{\{x:|x|\ge R\}}|\nabla u|^2dx\le \bigg(\frac{R}{2}\bigg)^{-2\alpha_0}
\int_{\{x:|x|\ge 1\}}|\nabla u|^2dx.
\end{equation}

\medskip


{\bf Step 3}. In this step, as opposed to the previous step, we deal with the issue that we do not know a priori that
$\eta^2(u-\bar u)$ is an admissible test function for (\ref{L3.13}). This allows us to rigorously deduce the continuum version
(\ref{L3.14b}) for $R\ge 1$, and where $u$ is now specified to be a partial derivative of the Green function,
i.e.\ $u(x)=\nabla_{i,x'} G(x,x')|_{x'=0}$ with $i=1,\cdots,d$. 
More precisely, we worry about the \emph{decay} at 
$|x|\to\infty$ --- we do not worry about local smoothness since anyway, we will
apply the argument to the discrete case in the next step. As in Step 1, we circumvent the problem of decay
through approximation by the \emph{periodic} problem. Indeed, we consider
the \emph{periodic} continuum elliptic Green function $G_L(x,x')=G_L(a,x,x')$ of period $L$. 
Up to additive constants, it is characterized by the weak equation
\begin{equation}\label{L3.16}
\int_{[-\frac{L}{2},\frac{L}{2})^d}\nabla\zeta(x)\cdot a(x)\nabla_xG_L(x,x')dx
=\zeta(x')-L^{-d}\int_{[-\frac{L}{2},\frac{L}{2})^d}\zeta(x)dx
\end{equation}
for all periodic $\zeta(x)$.
We note that $u_L(x)=\nabla_{i,x'} G_L(x,x')|_{x'=0}$ thus is characterized by
\begin{equation}\label{L3.12}
\int_{[-\frac{L}{2},\frac{L}{2})^d}\nabla\zeta(x)\cdot a(x)\nabla_x u_L(x)dx
=\nabla_i\zeta(0).
\end{equation}
Since $\nabla\nabla G_L$ converges in the sense of distributions to $\nabla\nabla G$ as $L\to\infty$, 
it is enough to show (\ref{L3.12}) implies
\begin{equation}\label{L3.17}
\int_{[-\frac{L}{2},\frac{L}{2})^d\cap\{x:|x|\ge R\}}|\nabla u_L|^2dx
\le C(d,\lambda)R^{-2\alpha_0}\int_{[-\frac{L}{2},\frac{L}{2})^d\cap\{x:|x|\ge 1\}}|\nabla u_L|^2 dx
\end{equation}
for $1\le R\le C(d)L$. 
Indeed we can estimate the right-hand side of \eqref{L3.17} using \eqref{energy_G_L} and 
apply weak lower semi-continuity to take the limit as $L\to\infty$ on the left-hand side to obtain \eqref{L3.14b}.
Now, disregarding smoothness issues, 
$\eta^2(u_L-\overbar{u_L})$ is an admissible test function for
(\ref{L3.12}). The argument for (\ref{L3.17}) is identical to the one in Step 2.

\medskip


{\bf Step 4}. Rigorous derivation of (\ref{L3.14}) for $R\ge C(d)$. In this step, we indicate
the modifications in Step 2 (or rather Step 3) that are necessary to treat the discrete case.
The first modification results from the fact that Leibniz rule and thus the neat identity
(\ref{L3.20}) does not hold anymore. However, we claim that the estimate (\ref{L3.21})
survives in form of
\begin{equation}\label{L3.23}
\nabla(\eta^2 (u-\bar u))(e) a(e)\nabla u(e)\ge \lambda(\nabla(\eta (u-\bar u))(e))^2
-([u](e)-\bar u)^2(\nabla\eta(e))^2,
\end{equation}
where we denote by $[u]([x,x+e_i])=\frac{1}{2}(u(x)+u(x+e_i))$ the local average of $u$ along
each edge $e=[x,x+e_i]$.
Indeed, since $\lambda\le a(e)\le 1$ is elliptic, this follows from the simple inequality on 4 \emph{numbers}
$\eta=\eta(x)$, $\tilde\eta=\eta(x+e_i)$, $v=u(x) - \bar u$, and $\tilde v=u(x+e_i) - \bar u$:
\begin{multline*}
(\eta^2v-\tilde\eta^2\tilde v)(v-\tilde v)-(\eta v-\tilde\eta\tilde v)^2
=-(\eta-\tilde\eta)^2v\tilde v\ge-(\eta-\tilde\eta)^2({\textstyle\frac{1}{2}}(v+\tilde v))^2.
\end{multline*}
Hence, if $\eta(x)$ denotes the (slightly narrower) cut-off function for $\{x:|x|\ge 2R-2\}$ 
in $\{x:|x|\ge R+2\}$ (which is possible for $R\ge 5$), from (\ref{L3.23}) we obtain the following
substitute of (\ref{L3.22})
\begin{align}\nonumber
\sum_{e:|e-b|\ge 2R}|\nabla u(e)|^2
&\le C(d,\lambda)R^{-2}\sum_{e:R+1\le |e-b|\le 2R-1}
([u](e)-\bar u)^2\\
\label{caccioppoli}
&\le C(d,\lambda)R^{-2}\sum_{x:R\le |b-x|\le 2R}
(u(x)-\bar u)^2.
\end{align}

\medskip

The second modification comes from the fact that we need a {\it discrete} version
of the Poincar\'e estimate with mean value zero on the annulus $\mathbb{Z}^d\cap\{R\le|x|\le 2R\}$,
which obviously holds with a constant $C(d)R^{2}$ provided that $R\ge C(d)$. 



\subsection{Proof of Lemma \ref{L2}}\label{ssec:L2}

{\bf Step 1}. In this first step, we consider two coefficient fields $\tilde a, a\in \Omega$ and their associated
Green functions $\tilde G = G(\tilde a; \cdot,\cdot)$ and $G = G(a; \cdot,\cdot)$, respectively. We claim that if $\tilde a$ and $a$ differ only at some edge $e\in\E^d$, then we
have that:
\begin{align}
 \tilde G(x,x') - G(x,x') &= (a(e) - \tilde a(e)) \nabla \tilde G(x,e) \nabla G(e,x'),\label{L2.26,5}\\
 \nabla \tilde G(b,x') - \nabla G(b,x') &= (a(e) - \tilde a(e)) \nabla \nabla \tilde G(b,e) \nabla G(e,x'),\label{L2.27}\\
 \nabla \nabla \tilde G(b,b') - \nabla \nabla G(b,b') &= (a(e) - \tilde a(e)) \nabla \nabla \tilde G(b,e) \nabla \nabla G(e,b').\label{L2.1}
\end{align}
Indeed, the difference satisfies the equation
\[
 \nabla^* \tilde a \nabla(\tilde G - G)(\cdot,x') = \nabla^* (a - \tilde a) \nabla G(\cdot, x')
\]
Since by assumption $\tilde a(b) = a(b)$ for all edges $b\neq e$, the Green function representation \eqref{Green_dist} immediately yields
\eqref{L2.26,5}. Differentiating \eqref{L2.26,5} then yields \eqref{L2.27} and \eqref{L2.1}.

\medskip

{\bf Step 2}. In this step, we derive the following estimate on the oscillations:
\begin{align}
 \osc_{a(e)} G(x,x') &\le 4\big({\textstyle 1 + \frac{1}{\lambda}}\big) |\nabla G(x,e)||\nabla G(e,x')|,\label{L2:osc1}\\
 \osc_{a(e)} \nabla G(b,x') &\le 4\big({\textstyle 1 + \frac{1}{\lambda}}\big) |\nabla\nabla G(b,e)| |\nabla G(e,x)|,\label{L2:osc2}\\
 \osc_{a(e)} \nabla\nabla G(b,b') &\le 4\big({\textstyle 1 + \frac{1}{\lambda}}\big) |\nabla\nabla G(b,e)||\nabla\nabla G(e,b')|.
\label{L2:osc3}
\end{align}
To do so, we first show that for any edge $e$, the dependence of $\nabla G(e,\cdot)$ on the value of $a(e)$ of
the conductivity is mild in the sense that
\begin{align}
 |\nabla \tilde G(e,x') - \nabla G(e,x')| &\le {\frac{1}{\lambda}} |\nabla G(e,x')|,\label{L2:Gronwall}\\
 |\nabla \nabla \tilde G(e,b') - \nabla\nabla G(e,b')| &\le {\frac{1}{\lambda}} |\nabla\nabla G(e,b')|,\label{L2:Gronwall2}
\end{align}
where $\tilde G$ and $G$ are given in Step~1.
This indeed follows from letting $b=e$ in \eqref{L2.27} and \eqref{L2.1} and recalling the a priori estimate $|\nabla\nabla \tilde G(e,e)|
\le \lambda^{-1}$ from  \eqref{L3.19}.
We turn to the proof of \eqref{L2:osc2}. It is clear that for any $a\in\Omega$, there exist $\tilde a_1,\tilde a_2\in \Omega$ with $\tilde
a_1(b) = a(b) = \tilde a_2(b)$ for all $b\neq e$ and associated Green functions $\tilde G_1$ and $\tilde G_2$ such that
\begin{align*}
 \osc_{a(e)} G(x,x') &\le 2 |\tilde G_1(x,x') - \tilde G_2(x,x')|\\
 &\le 2 |\tilde G_1(x,x') - G(x,x')| + 2 |G(x,x') - \tilde G_2(x,x')|.
\end{align*}
We insert \eqref{L2.26,5} with $\tilde a := \tilde a_i$, $i=1,2$, into this estimate to obtain that
\[
 \osc_{a(e)} G(x,x') \le 2 |\nabla \tilde G_1(x,e)| |\nabla G(e,x')| + 2 |\nabla \tilde G_2(x,e)| |\nabla G(e,x')|
\]
Consequently,  symmetry $\nabla \tilde G_i(x,e)=\nabla \tilde G_i(e,x)$ and estimate \eqref{L2:Gronwall} yield
\[
 \osc_{a(e)} G(x,x') \le 4\big({\textstyle 1 + \frac{1}{\lambda}}\big) |\nabla G(x,e)| |\nabla G(e,x')|.
\]
This proves \eqref{L2:osc1}. The estimates \eqref{L2:osc2} and \eqref{L2:osc3} follow similarly using \eqref{L2:Gronwall2}.

\medskip

{\bf Step 3}. In this step, we rephrase Lemma \ref{L3}, more precisely \eqref{L2.17}, 
in a way more suitable for its application in Step 4. More specifically, we claim
that there exists a weight exponent 
$\alpha(d,\lambda)>0$ such that
\begin{equation}\label{L2.16}
 \sup_{a\in\Omega}\sum_{e}\big|(|e-b|+1)^\alpha\nabla\nabla G(e,b)\big|^{2q}\le C(d,\lambda,q),
\end{equation}
for all $q \ge 1$. In fact, we claim that
\begin{equation}\label{L2.19}
\alpha:=\frac{1}{2}\alpha_0
\end{equation}
does the job. 
Because of $q\ge 1$, and thus $\ell^{2}(\E^d)\subset\ell^{2q}(\E^d)$, we have
\[
\sum_{e}\big|(|e-b|+1)^\alpha\nabla\nabla G(e,b)\big|^{2q}\le\Big(\sum_{e}\big|(|e-b|+1)^\alpha\nabla\nabla G(e,b)\big|^{2}\Big)^q.
\]
Using a dyadic decomposition, we see
\begin{align*}
&{\sum_{e}\big|(|e-b|+1)^\alpha\nabla\nabla G(e,b)\big|^{2}}\\
&=|\nabla\nabla G(b,b)|^2+\sum_{n=0}^\infty\sum_{e:2^{n-1}\le|e-b|<2^n} \big|(|e-b|+1)^\alpha\nabla\nabla G(e,b)\big|^{2}\\
&\le|\nabla\nabla G(b,b)|^2+\sum_{n=0}^\infty 2^{2\alpha (n+1)}\sum_{e:2^{n-1}\le|e|<2^n}|\nabla\nabla G(e,b)|^{2}.
\end{align*}
We now may appeal to \eqref{L2.17}  to obtain
\begin{equation}
 \sum_{e}\big|(|e-b|+1)^\alpha\nabla\nabla G(e,b)\big|^{2}\le C(d,\lambda)\Big(1+\sum_{n=0}^\infty 2^{2\alpha (n+1)} 2^{-2\alpha_0 n}\Big)
\stackrel{(\ref{L2.19})}{\le} C(d,\lambda).
\end{equation}

\medskip

{\bf Step 4}. In this step, we establish the first statement of Lemma \ref{L2},
namely (\ref{L2.25}). More precisely, we claim that for $p\ge \max\{\frac{d}{\alpha},1\}$ with $\alpha$ chosen in Step 3 and all $b,b'\in\E^d$:
\begin{multline}
 (|b-b'|+1)^{2pd} \Big\langle\bigg(\sum_{e}\Big(\osc_{a(e)} \nabla\nabla G(b,b') \Big)^2\bigg)^p\Big\rangle\\
 \le C(d,\lambda,p) \sup_{e, e'} \Big\{ (|e-e'|+1)^{2pd} \langle|\nabla\nabla G(e,e')|^{2p}\rangle\Big\}.\label{L2.15}
\end{multline}
Indeed, we first square \eqref{L2:osc3} and sum over $e$:
\[
\sum_e \Big(\osc_{a(e)} \nabla\nabla G(b,b') \Big)^2 \le C(\lambda) \sum_e |\nabla\nabla G(b,e)|^2|\nabla\nabla G(e,b')|^2.
\]
After taking the $p$-th power, we split the sum into its contributions over $\{e:|e-b|\le |e-b'|\}$ and $\{e:|e-b| > |e-b'|\}$ to obtain
\begin{multline}\label{step6_sumsplit}
 \bigg(\sum_e \Big(\osc_{a(e)} \nabla\nabla G(b,b') \Big)^2\bigg)^p\\
 \le C(\lambda,p) \Bigg(\bigg(\sum_{e:|e-b|\le|e-b'|}|\nabla\nabla G(b,e)|^2|\nabla\nabla G(e,b')|^2\bigg)^p\\
 + \bigg(\sum_{e:|e-b|\ge|e-b'|}|\nabla\nabla G(b,e)|^2|\nabla\nabla G(e,b')|^2\bigg)^p
\Bigg).
\end{multline}
We first bound the first term. To this end, we smuggle in a weight $(|e-b|+1)^{2\alpha}$ with $\alpha=\alpha(d,\lambda)$ from Step 3 and apply
H\"older's inequality
 with $p$ and its dual exponent $q$ (i.e.\ $\frac{1}{p} + \frac{1}{q} = 1$):
\begin{multline*}
 \bigg(\sum_{e:|e-b|\le|e-b'|}|\nabla\nabla G(b,e)|^2|\nabla\nabla G(e,b')|^2\bigg)^p\\
 \le \bigg(\sum_{e:|e-b|\le|e-b'|} \Big( (|e-b|+1)^\alpha |\nabla\nabla G(b,e)| \Big)^{2q} \bigg)^{p-1} \\
 \times \sum_{e:|e-b|\le|e-b'|} \Big( (|e-b|+1)^{-\alpha} |\nabla\nabla G(e,b)| \Big)^{2p}.
\end{multline*}
The first term on the right-hand side is bounded by Step 3, that is (\ref{L2.16}). After taking the expectation, we
smuggle in another weight $(|e-b'|+1)^{2pd}$ and take the supremum over appropriate terms to obtain
\begin{multline*}
 \Big\langle\sum_{e:|e-b|\le|e-b'|} \Big( (|e-b|+1)^{-\alpha} |\nabla\nabla G(e,b')| \Big)^{2p} \Big\rangle \\\le
 \bigg(\sum_{e:|e-b|\le|e-b'|} (|e-b|+1)^{-2p\alpha} (|e-b'|+1)^{-2pd}\bigg)\\
 \times \sup_{e'}\Big\{ (|e'-b'|+1)^{2pd} \big\langle|\nabla\nabla G(e',b')|^{2p} \big\rangle \Big\}.
\end{multline*}
Since $|e-b|\le|e-b'|$ implies $|e-b'|\ge\frac{1}{2}|b-b'|$, we find for the first r.-h.\ s.\ factor that
\begin{multline*}
 \sum_{e:|e-b|\le|e-b'|} (|e-b|+1)^{-2p\alpha} (|e-b'|+1)^{-2pd}\\
 \le \big({\textstyle\frac{1}{2}}|b-b'|+1\big)^{-2pd} \sum_{e} (|e-b|+1)^{-2p\alpha}.
\end{multline*}
Since by assumption $2p\alpha \ge 2d > d$, we obtain for the last factor
\[
\sum_{e\in\E^d}(|e-b|+1)^{-2p\alpha}\le C(d).
\]
Combining these estimates yields the bound
\begin{multline*}
 \Big\langle\Big(\sum_{e:|e-b|\le|e-b'|}|\nabla\nabla G(b,e)|^2|\nabla\nabla G(e,b')|^2\Big)^p\Big\rangle\\ \le
\bigg(C(d,\lambda,p) (|b-b'|+1)^{-d} \sup_{e,e'}\Big\{ (|e-e'|+1)^{d}
 \big\langle |\nabla\nabla G(e,e')|^{2p} \big\rangle^{\frac{1}{2p}} \Big\}\bigg)^{2p},
\end{multline*}
i.e.\ the expectation of the first term on the right-hand side of \eqref{step6_sumsplit} is bounded as desired.
 The second term in \eqref{step6_sumsplit} can be dealt with exactly as the first term by simply
exchanging the roles of $b$ and $b'$.

\medskip

{\bf Step 5}. Like in Step 3, we rephrase Lemma \ref{L3}, this time (\ref{I.2}),
in a way more suitable for its application in Step 6. We claim that for any
integrability exponent $q\ge 1$ and any weight exponent $\beta>0$ we have
\begin{equation}\label{L2.70}
 \sup_{a\in\Omega}\sum_{e}\big|(|e|+1)^{-\beta}\nabla G(e,0)\big|^{2q} \le C(d,\lambda,q,\beta)
\end{equation}
 We note that by (\ref{I.2}) we have as soon as $\beta>0$:
\begin{align}\label{L2.42}
 &\sum_{e}|(|e|+1)^{-\beta}\nabla G(e,0)|^{2q}\le\Big(\sum_{e}|(|e|+1)^{-\beta}\nabla G(e,0)|^{2}\Big)^q\nonumber\\
 &\le\Big(\sum_{i=1}^d |\nabla G([0,e_i],0)|^2
+\sum_{n=0}^\infty2^{-q\beta n}\sum_{e:2^n\le |e|<2^{n+1}}|\nabla G(e,0)|^2\Big)^q\\
 &\stackrel{(\ref{I.2})}{\le} C(d,\lambda,\beta).\nonumber
\end{align}

\medskip

{\bf Step 6}. In this step we establish the second conclusion of Lemma \ref{L2},
namely (\ref{L2.24}). More precisely, we show that for any integrability exponent
$p<\infty$ at least as large as in Step 3 and for any weight exponent $\beta>0$ sufficiently small such that
\begin{equation}\label{L2.41}
2p(\beta-d)+d<0
\end{equation}
we have
\begin{align}
 &(|b-x|+1)^{d-1} \Big\langle \bigg(\sum_e \Big(\osc_{a(e)}\nabla G(b,x) \Big)^2 \bigg)^p\Big\rangle^{\frac{1}{2p}}\nonumber\\
 &\le C(d,\lambda,p,\beta)\bigg(\sup_{e,x'} \Big\{(|e-x'|+1)^{d-1} \langle|\nabla G(e,x')|^{2p}\rangle^{\frac{1}{2p}} \Big\}\label{L2.40}\\
 &\qquad+(|b-x|+1)^{\beta-1+\frac{d}{2p}} \sup_{e,e'} \Big\{(|e-e'|+1)^{d}\langle|\nabla\nabla
G(e,e')|^{2p}\rangle^{\frac{1}{2p}}\Big\}\bigg),
\nonumber
\end{align}
for all $x\in\Z^d$ and $b\in\E^d$, where $C(d,\lambda,p,\beta)$ denotes a generic constant that only depends on $d$, $\lambda$,
$p$, and $\beta$. We note that by choosing $\beta$ small and $p$ large,
the exponent $\beta-1+\frac{d}{2p}$ can be made to be non-positive (in fact, as close to $-1$ as we want), which
proves (\ref{L2.24}).
In order to establish (\ref{L2.40}), we first square (\ref{L2:osc2}) and sum over $e\in\E^d$ to obtain that
\[
\sum_e \Big(\osc_{a(e)}\nabla G(b,x) \Big)^2 \le C(\lambda) \sum_e|\nabla\nabla G(b,e)|^2|\nabla G(e,x)|^2.
\]
We now split the sum over $e$:
\begin{align*}
&\sum_e|\nabla\nabla G(b,e)|^2|\nabla G(e,x)|^2\nonumber\\
&\le C(d,\lambda)
\Big(\sum_{e:|e-x|\ge\frac{1}{2}|b-x|}+\sum_{e:|e-x|<\frac{1}{2}|b-x|}\Big)|\nabla\nabla G(b,e)|^2|\nabla G(e,x)|^2.
\nonumber
\end{align*}
Since $|e-x|<\frac{1}{2}|b-x|$ implies $|e-b| > \frac{1}{2}|b-x|$, it follows
\begin{align}
 &\sum_e|\nabla\nabla G(b,e)|^2|\nabla G(e,x)|^2\nonumber\\
 &\le C(d,\lambda)\Big( \sum_{e:|e-x|\ge\frac{1}{2}|b-x|}|\nabla\nabla G(b,e)|^2|\nabla G(e,x)|^2\nonumber\\
 &\qquad\qquad\qquad+\sum_{e:|e-b| > \frac{1}{2}|b-x|}|\nabla G(e,x)|^2|\nabla\nabla G(b,e)|^2\Big).\label{L2.20}
\end{align}
We start by treating the first term on the r.-h.\ s. of (\ref{L2.20}) in an analogous way
to Step 4. For that purpose, let $\alpha$ be as in Step 3.
We smuggle in the weight $(|e-b|+1)^\alpha$ and apply H\"older's inequality with $p$ and $q$ such that $\frac{1}{p} + \frac{1}{q} = 1$:
\begin{multline*}
 \bigg(\sum_{e:|e-x|\ge\frac{1}{2}|b-x|}|\nabla\nabla G(b,e)|^2|\nabla G(e,x)|^2\bigg)^p\\
 \le \bigg(\sum_{e}\big|(|e-b|+1)^\alpha\nabla\nabla G(b,e)\big|^{2q}\bigg)^{p-1}\\
 \times\sum_{e:|e-x|\ge\frac{1}{2}|b-x|}\big|(|e-b|+1)^{-\alpha}\nabla G(e,x)\big|^{2p}.
\end{multline*}
The first term was bounded by a constant $C(d,\lambda,p)$ in Step 3.
Now we take the expectation $\langle\cdot\rangle$ w.\ r.\ t.\ $a$ and then smuggle in a weight $(|e-x|+1)^{2p(d-1)}$ to
obtain as desired:
\begin{align}\nonumber
&\Big\langle\Big(\sum_{e:|e-x|\ge\frac{1}{2}|b-x|} |\nabla\nabla G(b,e)|^2|\nabla G(e,x)|^2\Big)^p\Big\rangle\nonumber\\
&\le \bigg(\sum_{e:|e-x|\ge\frac{1}{2}|b-x|}(|e-b|+1)^{-2p\alpha} (|e-x|+1)^{-2p(d-1)} \nonumber\\
&\qquad\qquad\times\sup_{e'} \Big\{ (|e'-x|+1)^{2p(d-1)} 
\langle|\nabla G(e',x)|^{2p}\rangle\Big\}\bigg)\nonumber\\
&\stackrel{(\ref{L2.16})}{\le} C(d,\lambda,p) \big(|b-x|+1\big)^{-2p(d-1)} \sup_{e'} \Big\{
(|e'-x|+1)^{2p(d-1)} \langle|\nabla G(e',x)|^{2p}\rangle\Big\},\label{L2.46}
\end{align}
where we have used that $2p\alpha > d$.

\medskip

We now address the second term on the r.-h.\ s. of (\ref{L2.20}) in a similar way,
just exchanging the roles of $\nabla G$ and $\nabla\nabla G$, of $b$ and $x$,
and of $\alpha$ and $-\beta$,
where the weight exponent $\beta>0$ needs to satisfy (\ref{L2.41}).
By H\"older's inequality we obtain:
\begin{align*}
 &\bigg(\sum_{e:|e-b|\ge\frac{1}{2}|b-x|}|\nabla G(e,x)|^2|\nabla\nabla G(b,e)|^2\bigg)^p\\
 &\le \bigg(\sum_{e}\big|(|e-x|+1)^{-\beta}\nabla G(e,x)\big|^{2q}\bigg)^{p-1}\\
 &\qquad\times\sum_{e:|e-b|\ge\frac{1}{2}|b-x|}|(|e-x|+1)^{\beta}\nabla\nabla G(b,e)|^{2p}.
\end{align*}
The first term is bounded by Step 5 in form of (\ref{L2.70}). 
Taking the expectation and
smuggling in a weight $(|e-b|+1)^{2pd}$ yields
\begin{align*}
 &\Big\langle \sum_{e:|e-b|\ge\frac{1}{2}|b-x|}|(|e-x|+1)^{\beta}\nabla\nabla G(b,e)|^{2p}\Big\rangle\\
 &\le\sum_{e:|e-b|\ge\frac{1}{2}|b-x|}(|e-x|+1)^{2p\beta}(|e-b|+1)^{-2pd}\\
 &\qquad\times\sup_{e'}\Big\{(|e'-b|+1)^{2pd}\langle|\nabla\nabla G(b,e')|^{2p}\rangle\Big\}.
\end{align*}
We note that by the triangle inequality in form of $|e-x|\le |e-b|+|b-x|$, in the range (\ref{L2.41}) the remaining sum is bounded as
follows:
\begin{align*}
&\sum_{e:|e-b|\ge\frac{1}{2}|b-x|}(|e-x|+1)^{2p\beta}(|e-b|+1)^{-2pd}\\
&\le C(p,\beta) \bigg((|b-x|+1)^{2p\beta}\sum_{e:|e-b|\ge\frac{1}{2}|b-x|}(|e-b|+1)^{-2pd}\\
&\qquad\qquad\qquad\qquad+\sum_{e:|e-b|\ge\frac{1}{2}|b-x|}(|e-b|+1)^{2p(\beta-d)}\bigg)\\
&\le C(d,p,\beta) (|b-x|+1)^{2p(\beta-d)+d}.
\end{align*}
Hence we have obtained
\begin{multline}\label{L2.45}
\Big\langle\bigg(\sum_{e:|e-b|\ge\frac{1}{2}|b-x|}
|\nabla G(e,x)|^2|\nabla\nabla G(b,e)|^2\bigg)^p\Big\rangle\\
\le C(d,\lambda,p,\beta)(|b-x|+1)^{2p(\beta-d)+d} \sup_{e,e'}\Big\{(|e-e'|+1)^{2pd}\langle|\nabla\nabla
G(e,e')|^{2p}\rangle\Big\}.
\end{multline}
In view of (\ref{L2.20}), the combination of (\ref{L2.46}) and (\ref{L2.45}) 
as well as taking the $2p$-th root yields (\ref{L2.40}).



\subsection{Proof of Theorem \ref{T}}\label{ssec:T}

We start with the proof of (\ref{T.1}). For this purpose, we fix $b,b'\in\E^d$ and $p<\infty$;
by Jensen's inequality, we may assume that $p\ge p_0$ with $p_0$ from Lemma \ref{L2}.
Applying Lemma \ref{L1} to $\zeta(a)=\nabla\nabla G(a;b,b')$ and inserting the estimate (\ref{L2.25}) of Lemma \ref{L2} yields (after
redefining $\delta$)
\begin{multline*}
 (|b-b'|+1)^d \langle|\nabla\nabla G(b,b')|^{2p}\rangle^\frac{1}{2p}\le C(d,\lambda,\rho,p,\delta) (|b-b'|+1)^d
\langle|\nabla\nabla G(b,b')|\rangle \\
+\delta\sup_{e,e'} \Big\{ (|e-e'|+1)^d \langle|\nabla\nabla G(e,e')|^{2p}\rangle^\frac{1}{2p}\Big\}.
\end{multline*}
We now insert (\ref{T.2}) and take the supremum over $b$ and $b'$:
\begin{multline*}
 \sup_{b,b'} \Big\{ (|b-b'|+1)^d \langle|\nabla\nabla G(b,b')|^{2p}\rangle^\frac{1}{2p}\Big\}\\ \le C(d,\lambda,\rho,p,\delta) 
+\delta\sup_{e,e'} \Big( (|e-e'|+1)^d \langle|\nabla\nabla G(e,e')|^{2p}\rangle^\frac{1}{2p}\Big).
\end{multline*}
Choosing $\delta=1/2$, we obtain (\ref{T.1}).
We deal with the objection that $\sup_{e,b}\{(|e-b|+1)^d\langle|\nabla\nabla G(e,b)|^{2p}\rangle^{1/(2p)}\}$
may be infinite by first working with the periodic Green function $G_L$ as
in the proof of Lemma \ref{L3} and then letting $L\to\infty$.

\medskip

We now turn to the proof of (\ref{T.3}). With help of the just established
(\ref{T.1}), we may upgrade the result of Lemma \ref{L2}, cf.\ (\ref{L2.24}), to 
\begin{multline}
(|b-x|+1)^{d-1} \Big\langle\bigg(\sum_{e} \Big(\osc_{a(e)} \nabla G(b,x) \Big)^2 \bigg)^p\Big\rangle^\frac{1}{2p}\\
\le C(d,\lambda,\rho,p)\bigg( \sup_{e,x'}\Big\{(|e-x'|+1)^{d-1} \langle|\nabla G(e,x')|^{2p}\rangle^\frac{1}{2p} \Big\}
+1\bigg).\label{T.4}
\end{multline}
We apply Lemma \ref{L1} to $\zeta=\nabla G(b,x)$ and insert \eqref{T.4} (after redefining $\delta$):
\begin{align*}
& (|b-x|+1)^{d-1} \langle|\nabla G(b,x)|^{2p}\rangle^\frac{1}{2p}\\
&\le C(d,\lambda,\rho,p,\delta) (|b-x|+1)^{d-1}\langle|\nabla G(b,x)|\rangle\\
&\qquad\qquad+\delta\bigg(\sup_{e,x'}\Big\{ (|e-x'|+1)^{d-1} \langle|\nabla G(e,x')|^{2p}\rangle^\frac{1}{2p} \Big\} + 1
\bigg).
\end{align*}
We now insert (\ref{I.1}) and take the supremum over $b$ and $x$:
\begin{multline*}
\sup_{b,x}\Big( (|b-x|+1)^{d-1} \langle|\nabla G(b,x)|^{2p}\rangle^\frac{1}{2p} \Big)\\
\le C(d,\lambda,\rho,p,\delta)+\delta\sup_{e,x'}\Big\{ (|e-x'|+1)^{d-1} \langle|\nabla G(e,x')|^{2p}\rangle^\frac{1}{2p}
\Big\}.
\end{multline*}
As before, letting $\delta=1/2$ yields \eqref{T.3}.



\subsection{Proof of Corollary \ref{C}}\label{ssec:C}

It is well known that an LSI implies a corresponding SG, see for instance \cite[Theorem 4.9]{GuionnetZegarlinski}. 
Indeed, using $\zeta^2=1+\epsilon f$ for some $f(a)$
in (\ref{LSI}) and expanding to second order in $\epsilon\ll 1$ one obtains
\[
\langle (f - \langle f \rangle)^2 \rangle\le\frac{1}{\rho}\Big\langle\sum_e \Big(\osc_{a(e)}f\Big)^2 \Big\rangle.
\]
As in Step 2 of the proof of Lemma \ref{L1}, see also \cite[Lemma 11]{GloriaNeukammOtto}, it follows that
\begin{equation}\label{SGp}
\langle |f - \langle f \rangle|^{2p} \rangle\le C(\rho,p) \Big\langle\bigg(\sum_e \Big(\osc_{a(e)}f\Big)^2 \bigg)^p\Big\rangle.
\end{equation}
We fix $x\in\mathbb{Z}^d$ and apply this inequality to $f(a)=G(a;x,0)$ and use \eqref{L2:osc1} from the proof
of Lemma \ref{L2}, i.e.\
\begin{equation}\nonumber
 \osc_{a(e)} G(x,0)\le C(\lambda)|\nabla G(x,e)||\nabla G(e,0)|,
\end{equation}
to obtain
\begin{equation}\nonumber
\big\langle \big| G(x,0) -\langle G(x,0) \rangle \big|^{2p} \big\rangle^{\frac{1}{p}}\le
C(\lambda,\rho,p)\Big\langle\bigg(\sum_e |\nabla G(x,e)|^2|\nabla G(e,0)|^2 \bigg)^p\Big\rangle^{\frac{1}{p}}.
\end{equation}
The triangle inequality in $\langle (\cdot)^p\rangle^{1/p}$ yields
\[
 \big\langle \big| G(x,0) -\langle G(x,0) \rangle \big|^{2p} \big\rangle^{\frac{1}{p}}\le
C(\lambda,\rho,p) \sum_e \Big\langle |\nabla G(x,e)|^{2p}|\nabla G(e,0)|^{2p} \Big\rangle^{\frac{1}{p}}.
\]
Using the Cauchy-Schwarz inequality in $\langle\cdot\rangle$ and appealing to
stationarity, we obtain
\begin{align*}
& \big\langle \big| G(x,0) -\langle G(x,0) \rangle \big|^{2p} \big\rangle^{\frac{1}{p}}\\
&\le C(\lambda,\rho,p)\sum_e \langle|\nabla G(x,e)|^{4p}\rangle^\frac{1}{2p}\langle|\nabla G(e,0)|^{4p}\rangle^\frac{1}{2p}\\
&=C(\lambda,\rho,p)\sum_e \langle|\nabla G(e-x,0)|^{4p}\rangle^\frac{1}{2p}\langle|\nabla G(e,0)|^{4p}\rangle^\frac{1}{2p},
\end{align*}
where we recall that $e-x\in\E^d$ is the edge $e$ shifted by $x$ and $\nabla$ always falls on the edge variable.
Into this estimate, we insert the result of Theorem \ref{T}:
\begin{equation}\label{C.2}
\big\langle \big| G(x,0) -\langle G(x,0) \rangle \big|^{2p} \big\rangle^{\frac{1}{p}}\le
C(d,\lambda,\rho,p)\sum_e \big( (|e-x|+1)(|e|+1) \big)^{2(1-d)}.
\end{equation}

\medskip

We now turn to the sum on the r.-h.\ s. of (\ref{C.2}): By symmetry, we have
\begin{equation}\label{C.4}
\sum_e \big( (|e-x|+1)(|e|+1) \big)^{2(1-d)}\le 2\sum_{e:|e-x|\le |e|}\big( (|e-x|+1)(|e|+1) \big)^{2(1-d)}.
\end{equation}
We note that in the case of $d>2$ we have $2(1-d)<-d$ so that
\begin{equation}\label{T.5}
\sum_e(|e|+1)^{2(1-d)}\le C(d)<\infty.
\end{equation}
Since $|e-x|\le|e|$ implies $|e|\ge\frac{1}{2}|x|$ we thus have as desired for (\ref{C.4})
\begin{align}
 &\sum_{e:|e-x|\le |e|}(|e-x|+1)^{2(1-d)}(|e|+1)^{2(1-d)}\nonumber\\
 &\le ({\textstyle\frac{1}{2}}|x|+1)^{2(1-d)}\sum_{e}(|e-x|+1)^{2(1-d)}\nonumber\\
 &\stackrel{(\ref{T.5})}{\le} C(d)(|x|+1)^{2(1-d)}.\label{C.5}
\end{align}
We now turn to the case of $d=2$. In this case, we split the sum on the r.~h.~s. of (\ref{C.4})
according to
\begin{align*}
\sum_{e:|e-x|\le |e|}
&=\sum_{e:|e-x|\le |e|\;\mbox{\scriptsize and}\;|e|\ge 2|x|}+\sum_{e:|e-x|\le |e|\;\mbox{\scriptsize and}\;|e|<
2|x|}\\
&\le\sum_{e:|e-x|\ge\frac{1}{2}|e|\;\mbox{\scriptsize and}\;|e|\ge 2|x|}
+\sum_{e:|e-x|\le 2|x|\;\mbox{\scriptsize and}\;|e|\ge \frac{1}{2}|x|},
\end{align*}
so that
\begin{align}
 &\sum_{e:|e-x|\le |e|}(|e-x|+1)^{-2}(|e|+1)^{-2}\nonumber\\
 &\le \sum_{e:|e|\ge 2|x|}({\textstyle\frac{1}{2}}|e|+1)^{-4}+
({\textstyle\frac{1}{2}}|x|+1)^{-2}\sum_{e:|e-x|\le 2|x|}(|e-x|+1)^{-2}\nonumber\\
 &\le C(|x|+1)^{-2}+C(|x|+1)^{-2}\log(|x|+2).\label{C.6}
\end{align}
Combining (\ref{C.5}) and (\ref{C.6}), we gather
\begin{align}
 &\sum_e(|e-x|+1)^{2(1-d)}(|e|+1)^{2(1-d)}\nonumber\\
 &\le C(d)(|x|+1)^{2(1-d)}\begin{cases}
			   1 &\text{for $d>2$}\\
			   \log(|x|+2) &\text{for $d=2$}
                          \end{cases}\Bigg\},\label{C:cases}
\end{align}
which we insert into (\ref{C.2}) to obtain (\ref{C.1}).



\subsection{Optimality of Corollary \ref{C} for $p=1$}\label{ssec:opt_C}

In this section we will show by formal calculations that Corollary \ref{C} is optimal by considering the regime $1 -
\lambda \ll 1$. Recall that the Green function satisfies $\nabla^* a \nabla G(\cdot,x') = \delta(\cdot-x')$. Now let
$a(e) = 1 + \epsilon \tilde a(e)$ for $\epsilon\ll1$, where $\tilde a$ is i.~i.~d.~with values at each edge taken in $[-1,1]$.
Furthermore we assume $\langle \tilde a(e) \rangle = 0$ as well as $\langle \tilde a(e)^2 \rangle = 1$. 
Note that this implies $a \in [1-\epsilon, 1+\epsilon] \subset [1/2, 3/2]$ (w.~l.~o.~g.~$\epsilon < 1/2$), 
but (by linearity of the equation in $a$) all results remain true with this new upper bound on $a$. Let us expand the
Green function corresponding to $a$ in powers of $\epsilon$:
\[
 G(x,y) = G_0(x,y) + \epsilon G_1(x,y) + \ldots\;.
\]
Substituting into the defining equation for $G$, we find that to zeroth order in $\epsilon$, we have
\[
 \nabla^* \nabla G_0(\cdot,x') = \delta(\cdot-x'),
\]
i.e.\ $G_0$ is the constant-coefficient Green function. Then to first order, it follows
\[
 \nabla^* \nabla G_1(\cdot,x') + \nabla^* \tilde a \nabla G_0(\cdot,x') = 0.
\]
Hence we have that
\[
 G_1(x,x') = -\sum_e \nabla G_0(x,e) \tilde a(e) \nabla G_0(e,x').
\]
Since $\langle \tilde a(e) \rangle = 0$, we deduce $\langle G_1 \rangle = 0$ and consequently
\begin{align*}
 \langle G_1(x,0)^2 \rangle &= \langle (G_1(x,0) - \langle G_1(x,0) \rangle)^2 \rangle\\
 &= \sum_{e,e'} \nabla G_0(x,e) \nabla G_0(x,e') \langle \tilde a(e) \tilde a(e') \rangle \nabla G_0(e,0) \nabla
G_0(e',0).
\end{align*}
Since the coefficients $\tilde a(x)$ are i.\ i.\ d.\ with variance 1, it follows
\[
 \langle (G_1(x,0) - \langle G_1(x,0) \rangle)^2 \rangle = \sum_{e} (\nabla G_0(x,e))^2 (\nabla G_0(e,0))^2.
\]
The behavior of the constant-coefficient Green function $G_0$ is well-known, cf.~\cite[Theorem 4.3.1]{LawlerLimic}, and yields
that $(\nabla G_0(e,0))^2$ scales like $(|e|+1)^{1-d}$ with a similar expression for $(\nabla G_0(x,e))^2$.  Hence we find that
\begin{align}\nonumber
\langle (G_1(x,0) - \langle G_1(x,0) \rangle)^2 \rangle &\le C(d) \sum_e \big((|e-x|+1) (|e|+1)\big)^{2(1-d)}
\;\;\text{and}\\
\langle (G_1(x,0) - \langle G_1(x,0) \rangle)^2 \rangle &\ge \frac{1}{C(d)} \sum_e \big((|e-x|+1) (|e|+1)\big)^{2(1-d)}.
\label{optim_scale_low}
\end{align}
Thus \eqref{C:cases} and \eqref{optim_scale_low} yield the upper bound
\[
 \langle (G_1(x,0) - \langle G_1(x,0) \rangle)^2 \rangle \le \frac{1}{C(d)} (|x|+1)^{2(1-d)}\begin{cases}
                                                             1 &\text{for $d>2$}\\
							     \log(|x|+2) &\text{for $d=2$}
                                                            \end{cases}\Bigg\}.
\]
If $d > 2$, a lower bound can be obtained by considering only the summand $e=[0,e_i]$ in
\eqref{optim_scale_low}. 
If $d=2$, we restrict the sum to all $e$ such that $|e| \le |x|$ and use $|e-x|\le2|x|$ in that region to obtain
\begin{align*}
 \langle (G_1(x,0) - \langle G_1(x,0) \rangle)^2 \rangle &\ge \frac{1}{C(d)} \sum_{e:|e|\le|x|} (|e-x|+1)^{-2}
(|e|+1)^{-2}\\
 &\ge \frac{1}{C(d)} (2|x|+1)^{-2} \sum_{e:|e|\le|x|} (|e|+1)^{-2}\\
 &\ge \frac{1}{C(d)} (|x|+1)^{-2} \log(|x|+2).
\end{align*}
Thus Corollary \ref{C} is indeed optimal in scaling.



\subsection{Proof of Corollaries \ref{Conlon_T1.2} and \ref{Conlon_T1.3}}\label{ssec:Conlon}

%
%
{\sc Proof of Corollary \ref{Conlon_T1.2}}

{\bf Step 1}. Proof in dimension $d > 2$. First of all, the triangle inequality in $\langle(\cdot)^r\rangle^{1/r}$
yields
\begin{equation}\label{Con_T1.2_eq1}
 \bigg\langle \bigg( \sum_x \big|u(x) - \langle u(x) \rangle \big|^p \bigg)^r
\bigg\rangle^\frac{1}{rp}
\le \bigg( \sum_x \Big\langle \big|u(x) - \langle u(x) \rangle \big|^{rp}
\Big\rangle^\frac{1}{r}
\bigg)^\frac{1}{p}.
\end{equation}
Since $u$ is the decaying solution of \eqref{u_eps} with compactly supported 
right-hand side $f$, it can be represented via the Green function:
\begin{equation}\label{Green_rep}
 u(x) = \sum_y G(x,y) f(y),
\end{equation}
Consequently, an application of the triangle inequality in $\langle(\cdot)^{rp}\rangle^{1/(rp)}$ yields
\begin{equation}\nonumber
 \Big\langle \big|u(x) - \langle u(x) \rangle \big|^{rp} \Big\rangle^\frac{1}{rp} \le \sum_y
\Big\langle \big| G(x,y) - \langle G(x,y) \rangle\big|^{rp} \Big\rangle^\frac{1}{rp}  |f(y)|,
\end{equation}
so that we may use Corollary \ref{C} to the effect of
\begin{equation}\label{Con_T1.2_eq2}
 \Big\langle \big|u(x) - \langle u(x) \rangle \big|^{rp} \Big\rangle^\frac{1}{rp} 
\le C(d,\lambda,\rho,r,p) \sum_y(|x-y|+1)^{1-d}|f(y)|.
\end{equation}
We now insert \eqref{Con_T1.2_eq2} in \eqref{Con_T1.2_eq1} to obtain
\begin{multline}\label{Con_T1.2_eq3}
 \bigg\langle \bigg( \sum_x \big|u(x) - \langle u(x) \rangle \big|^p \bigg)^r
\bigg\rangle^\frac{1}{rp}\\ \le C(d,\lambda,\rho,r,p) \Bigg( \sum_x \bigg( \sum_y (|x-y|+1)^{1-d}  |f(y)| \bigg)^p
\Bigg)^\frac{1}{p}.
\end{multline}
Now let us recall the Hardy-Littlewood-Sobolev inequality in $\R^d$, see \cite[Section 4.3]{LiebLoss} for a proof:
\[
 \Bigg( \int_{\R^d} \bigg( \int_{\R^d} |x-y|^{-\alpha} f(y) \;dy \bigg)^p \Bigg)^{\frac{1}{p}} \le C(d,\alpha,p) \bigg( \int_{\R^d} |f(y)|^q
\;dy \bigg)^{\frac{1}{q}}
\]
for all weight exponents $0<\alpha<d$ and
for all integrability exponents $1<p, q<\infty$ related by $1+\frac{1}{p} = \frac{\alpha}{d} + \frac{1}{q}$. 
A discrete version can easily be obtained by applying the continuum version to 
piecewise constant functions. We use the discrete version 
for $\alpha=d-1$, that is,
\begin{equation}\label{HLS}
  \Bigg( \sum_x \Big( \sum_y (|x-y|+1)^{1-d} |f(y)| \Big)^{p}
\Bigg)^{\frac{1}{p}}
 \le C(d,p) \Big( \sum_y |f(y)|^q \Big)^{\frac{1}{q}},
\end{equation}
in which case the relation turns as desired into $\frac{1}{p}+\frac{1}{d}=\frac{1}{q}$. Our assumption
$p\ge 2$ and $d>2$ ensure that $q$ is indeed admissible for Hardy-Littlewood-Sobolev in the sense of the
strict inequality $q>1$.

\medskip

{\bf Step 2}. Changes if $d=2$. In this case, using that $f(y)$ is supported in $\{y:|y|\le R\}$,
(\ref{Con_T1.2_eq3}) assumes the form
\begin{align*}
 &\bigg\langle \bigg( \sum_{x:|x|\le R} \big|u(x) - \langle u(x) \rangle \big|^p \bigg)^r\bigg\rangle^\frac{1}{rp}\\ 
 &\le C(d,\lambda,\rho,r,p) \Bigg( \sum_{x:|x|\le R} \bigg( \sum_{y:|y|\le R} (|x-y|+1)^{1-d} (\log^\frac{1}{2}|x-y|)  |f(y)|
\bigg)^p\Bigg)^\frac{1}{p}\\
 &\le C(d,\lambda,\rho,r,p) (\log^\frac{1}{2}R) \Bigg( \sum_{x} \bigg( \sum_{y} (|x-y|+1)^{1-d}  |f(y)| \bigg)^p\Bigg)^\frac{1}{p}.
\end{align*}
As in Step 1, it remains to apply the discrete
Hardy-Littlewood-Sobolev inequality, where we note that our assumption $p>2$ now ensures $q>1$
even for $d=2$.

\bigskip


{\sc Proof of Corollary \ref{Conlon_T1.3}}

{\bf Step 1}. In this step, we derive the estimate
\begin{multline}\label{Con_T1.3_step1}
 \Big\langle \Big| \sum_x (u(x)-\langle u(x)\rangle)g(x)\Big|^{r} \Big\rangle^{\frac{1}{r}}\\
 \le C(\rho,r) \Big\langle \bigg( \sum_e  \Big(\sum_x \sum_y \Big(\osc_{a(e)} G(x,y)\Big) |f(y)| |g(x)| \Big)^{2} \bigg)^{\frac{r}{2}}
\Big\rangle^{\frac{1}{r}}.
\end{multline}
%
Indeed, it follows from the representation \eqref{Green_rep} that
\begin{multline*}
 \Big\langle \Big| \sum_x \big( u(x) - \langle u(x) \rangle \big) g(x) \Big|^{r}
\Big\rangle^{\frac{1}{r}}\\
 = \Big\langle \Big| \sum_x \sum_y \big( G(x,y) - \langle G(x,y) \rangle \big)
 f(y) g(x) \Big|^{r} \Big\rangle^{\frac{1}{r}}.
\end{multline*}
Hence the $L^p$-version of SG \eqref{SGp}, with $2p$ replaced by $r$ (w.~l.~o.~g.~we may assume $r\ge 2$), yields
\begin{multline*}
 \Big\langle \Big| \sum_x \big( u(x) - \langle u(x) \rangle \big) g(x) \Big|^{r}
\Big\rangle^{\frac{1}{r}}\\
 \le C(\rho,r) \Big\langle \bigg( \sum_e  \Big(\osc_{a(e)}\sum_x \sum_y G(x,y) f(y) g(x) \Big)^{2} \bigg)^{\frac{r}{2}}
\Big\rangle^{\frac{1}{r}}.
\end{multline*}
Since the only dependence on the coefficients $a$ is through $G$, we may use sub-linearity of the oscillation to obtain
\eqref{Con_T1.3_step1}.

\medskip

{\bf Step 2.} In this step, we estimate the right-hand side of \eqref{Con_T1.3_step1} as follows:
\begin{align}\label{Con_T1.3_step2}
 &\Big\langle \bigg( \sum_e \Big( \sum_x \sum_y \Big(\osc_{a(e)} G(x,y)\Big) |f(y)| |g(x)| \Big)^{2} \bigg)^{\frac{r}{2}}
\Big\rangle^{\frac{1}{r}}\nonumber\\
 &\le C(d,\lambda,\rho,r) \Bigg( \sum_e \bigg( \sum_x (|e-x|+1)^{1-d} |g(x)| \bigg)^2\\
 &\qquad\qquad\qquad\qquad\qquad\times \bigg( \sum_y (|e-y|+1)^{1-d} |f(y)| \bigg)^2 \Bigg)^{\frac{1}{2}}.\nonumber
\end{align}
Indeed, expanding the square on the l.-h.\ s.\ of \eqref{Con_T1.3_step2} and inserting \eqref{L2:osc1}
yields
\begin{align*}
 &\Big( \sum_x \sum_y \Big(\osc_{a(e)} G(x,y)\Big) |f(y)| |g(x)| \Big)^{2}\\
 &\le C(\lambda) \sum_{x,x',y,y'} |\nabla G(x,e)| |\nabla G(x',e)| |\nabla G(e,y)| |\nabla G(e,y')|\\&\qquad\qquad\qquad\qquad\qquad\qquad
\times |g(x)| |g(x')| |f(y)| |f(y')|.
\end{align*}
Consequently we obtain by the triangle inequality w.\ r.\ t.\ $\langle |\cdot|^{r/2} \rangle^{2/r}$ that
\begin{align*}
 & \Big\langle \bigg( \sum_e \Big( \sum_x \sum_y \Big(\osc_{a(e)} G(x,y)\Big) |f(y)| |g(x)| \Big)^{2} \bigg)^{\frac{r}{2}}
\Big\rangle^{\frac{1}{r}}\\
 &\le C(d,\lambda) \Bigg(\sum_{e,x,x',y,y'} \langle |\nabla G(x,e)|^{\frac{r}{2}} |\nabla G(x',e)|^{\frac{r}{2}}
|\nabla G(e,y)|^{\frac{r}{2}} |\nabla G(e,y')|^{\frac{r}{2}}
\rangle^{\frac{1}{r}}\\&\qquad\qquad\qquad\qquad\qquad\qquad\qquad\times |g(x)| |g(x')|
|f(y)| |f(y')| \Bigg)^{\frac{1}{2}}.
\end{align*}
H\"older's inequality with respect to the ensemble $\langle \cdot \rangle$ and Theorem \ref{T} yield
\begin{align*}
 &\langle |\nabla G(x,e)|^{\frac{r}{2}} |\nabla G(x',e)|^{\frac{r}{2}} |\nabla G(e,y)|^{\frac{r}{2}} |\nabla G(e,y')|^{\frac{r}{2}}
\rangle^{\frac{2}{r}}\\
 &\le C(d,\lambda,\rho,r) \big((|e-x|+1) (|e-x'|+1) (|e-y|+1) (|e-y'|+1) \big)^{1-d}.
\end{align*}
Hence (\ref{Con_T1.3_step2}) follows from partly undoing the expansion of the square:
\begin{align*}
 &\Bigg( \sum_{e,x,x',y,y'} (|e-x|+1)^{1-d} (|e-x'|+1)^{1-d} (|e-y|+1)^{1-d}\\&\qquad\qquad\qquad\times
(|e-y'|+1)^{1-d} |g(x)| |g(x')| |f(y)| |f(y')| \Bigg)^{\frac{1}{2}}\\
 &= \Bigg( \sum_e \bigg( \sum_x (|e-x|+1)^{1-d} |g(x)| \bigg)^2
\bigg( \sum_y (|e-y|+1)^{1-d} |f(y)| \bigg)^2 \Bigg)^{\frac{1}{2}}.
\end{align*}

\medskip

{\bf Step 3.} Conclusion. 
An application of H\"older's inequality w.\ r.\ t.\ the sum over $e$ on the r.-h.\ s.\ of (\ref{Con_T1.3_step2}) yields
a bound by
\begin{multline}\label{Conlon_proof}
 \Bigg( \sum_e \bigg( \sum_x (|e-x|+1)^{1-d} |g(x)|
\bigg)^{\tilde p} \Bigg)^{\frac{1}{\tilde p}}\\ \times
 \Bigg( \sum_e \bigg( \sum_y (|e-y|+1)^{1-d} |f(y)|
\bigg)^{p} \Bigg)^{\frac{1}{p}},
\end{multline}
with $\tilde p$ and $p$ such that $\frac{2}{p} + \frac{2}{\tilde p} = 1$ to be chosen later.
We recall the Hardy-Littlewood-Sobolev inequality \eqref{HLS}, i.e.\
\[
 \Bigg( \sum_e \bigg( \sum_x (|e-x|+1)^{1-d} |f(x)| \bigg)^{p} \Bigg)^{\frac{1}{p}} \le C(d,q) \bigg( \sum_x
|f(x)|^q\bigg)^{\frac{1}{q}},
\]
if we choose $p$ such that $\frac{1}{q} = \frac{1}{d} + \frac{1}{p}$. (Here we require $1 < q < d$ so that in particular $q < \infty$.) The Hardy-Littlewood-Sobolev inequality
likewise yields
\[
 \Bigg(\sum_e \bigg( \sum_x (|e-x|+1)^{1-d} |g(x)| \bigg)^{\tilde p} \Bigg)^{\frac{1}{\tilde p}}
 \le C(d,q) \bigg( \sum_x |g(x)|^{\tilde q} \bigg)^{\frac{1}{\tilde q}},
\]
where $\frac{1}{\tilde q} = \frac{1}{d} + \frac{1}{\tilde p} = \frac{1}{d} + \frac{1}{2} - \frac{1}{p}$ and we require $1 < \tilde q < d$.
Inserting these
estimates into \eqref{Conlon_proof} and then into Steps 2 and 1 yields Corollary \ref{Conlon_T1.3}.



\subsection{Proof of Corollary \ref{C2}}\label{ssec:C2}

{\bf Step 1}. Let $u$ satisfy $\nabla^*a\nabla u = 0$ in $\{x:|x|\le R\}$. We claim that for any function $\eta(x)$ supported in $\{ x : |x| < R \}$, we obtain the representation
\begin{equation}\label{represent_eta_u}
 (\eta u)(x) = \sum_{e=[y,y']:|e| \le R} \big( u(y) \nabla G(e,x) a(e)  \nabla \eta(e) - G(y,x) \nabla \eta(e)
a(e) \nabla u(e) \big),
\end{equation}
where we sum over all edges in $\E^d$ of the form $[y,y']$ such that their midpoint is of distance at most $R$ from the origin.
We start by noting that even on the discrete level, some aspects of Leibniz rule survive, such as
\begin{align}\label{leibniz}
 &\nabla \zeta(e) a(e) \nabla (\eta u)(e) - \nabla (\eta \zeta)(e) a(e) \nabla u(e)\nonumber\\
 &= u(y) \nabla\zeta(e) a(e) \nabla \eta(e) - \zeta(y) \nabla\eta(e) a(e) \nabla u(e)
\end{align}
for any function $\zeta:\Z^d\to\R$ and $e=[y,y']\in\E^d$. Indeed, \eqref{leibniz} reduces to the elementary identity
\[
 (\tilde\zeta - \zeta)(\tilde\eta \tilde u - \eta u) - (\tilde \eta \tilde \zeta - \eta \zeta)(\tilde u - u) = u (\tilde\zeta -
\zeta)(\tilde\eta - \eta) - \zeta (\tilde\eta - \eta)(\tilde u - u).
\]
We integrate \eqref{leibniz}:
\begin{equation}\label{int_leibniz}
 \sum_e \nabla \zeta a \nabla (\eta u) - \sum_e \nabla (\eta \zeta) a \nabla u = \sum_e \big(u \nabla\zeta a \nabla \eta
- \zeta\nabla\eta a \nabla u\big)
\end{equation}
and use it for $\zeta = G(\cdot, x)$. By definition of $G$, the first term on the l.-h.\ s.\ of \eqref{int_leibniz} yields $(\eta u)(x)$. 
Since $\eta G(\cdot, x)$ is supported in $\{ y : |y| \le R \}$, the second term on the l.-h.\ s.\ of \eqref{int_leibniz} vanishes. This
completes the step.

\medskip

{\bf Step 2}. We now use the representation obtained in Step 1 to obtain bounds on the gradient of $u$ and consequently on the
$\alpha$-H\"older norm of $u$. Specifically, we claim that
\begin{multline}\label{C2_de_giorgi_green}
 \Bigg(\frac{\sup_{x:|x|\le \frac{R}{8}} \frac{|u(x) - u(0)|}{|x|^\alpha}}{\frac{1}{R^\alpha} \sup_{x:|x|\le R} |u(x)|} \Bigg)^p\\
 \le C(d,\lambda,p) R^{\alpha p} R^{-p} \Big( R^{d(p-1)} \sum_{e:|e|\le \frac{R}{8}} \sum_{b:\frac{R}{4} \le |b|\le \frac{R}{2}}
|\nabla\nabla G(e,b)|^p\\
 + R^{d(p-1)-p} \sum_{e:|e|\le \frac{R}{8}} \sum_{x:\frac{R}{4} \le |x|\le \frac{R}{2}} |\nabla G(e,x)|^p \Big),
\end{multline}
if $\alpha<1$ and $p>d$ are related by $\alpha p = p-d$.
To this end, we choose a cut-off function $\eta$ for $\{ x : |x| \le \frac{R}{4} + 1 \}$ in $\{ x : |x| \le \frac{R}{2} - 1 \}$ (w.~l.~o.~g.~$R > 8$). 
We restrict to $|x| \le \frac{R}{4}$ and take the derivative of \eqref{represent_eta_u} along the edge $b$ to obtain
\[
 \nabla u(b) = \sum_{e=[y,y']\in\E^d} \big( u(y) \nabla\nabla G(e,b)  a(e)  \nabla \eta(e) - (\nabla \eta(e) a(e) \nabla u(e))
\nabla G(y,b) \big).
\]
This implies
\begin{equation}\label{C2_grad_estimate}
 |\nabla u(b)| \le \frac{C(d)}{R} \sum_{e=[y,y']:\frac{R}{4} \le |y|,|y'| \le \frac{R}{2}} \big( |u(y)| |\nabla\nabla G(e,b)| +
|\nabla G(y,b)| |\nabla u(e)| \big).
\end{equation}
Applying H\"older's inequality and summing the $p$-th power of \eqref{C2_grad_estimate}, we obtain
\begin{multline}\label{C2_grad_estimate2}
 \sum_{b:|b|\le \frac{R}{8}}|\nabla u(b)|^p\\ \le C(d,p) R^{-p} \bigg( \Big(\sum_{y:|y|\le\frac{R}{2}} |u(y)|^q\Big)^{p-1}
\sum_{b:|b|\le\frac{R}{8}}\sum_{e:\frac{R}{4} \le |e|\le\frac{R}{2}}|\nabla\nabla G(e,b)|^p\\
 + \Big(\sum_{e:|e|\le\frac{R}{2}}|\nabla u(e)|^q\Big)^{p-1} \sum_{b:|b|\le\frac{R}{8}}\sum_{y:\frac{R}{4} \le |y|\le\frac{R}{2}} |\nabla
G(y,b)|^p \bigg),
\end{multline}
where $q$ is the dual H\"older exponent of $p$. 
Now we apply the following (discrete) Sobolev inequality: If $\alpha < 1$ and $p > d$ are related by
\begin{equation}\label{C2_sobolev_index}
 \alpha = 1 - \frac{d}{p},
\end{equation}
then we have that
\begin{equation}\label{morrey}
 \sup_{x:|x|\le \frac{R}{8}}  \frac{|u(x) - u(0)|}{|x|^\alpha} \le C(d,p) \bigg( \sum_{b:|b|\le \frac{R}{8}} |\nabla u(b)|^p
\bigg)^{\frac{1}{p}}.
\end{equation}
(This discrete version can easily be derived from its continuum version by extending $u$ to a piecewise linear function on a 
triangulation subordinate to the lattice.) Therefore the left-hand side of \eqref{C2_grad_estimate2} bounds the
$\alpha$-H\"older norm as desired, albeit over a smaller ball.

\medskip

Let us now turn to the right-hand side of \eqref{C2_grad_estimate2}. We trivially have that
\begin{equation}\label{C2_Lq_est}
 \bigg( \sum_{y:|y| \le \frac{R}{2}} |u(y)|^q \bigg)^{p-1} \le C(d,p) R^{d(p-1)} \Big( \sup_{x:|x|\le \frac{R}{2}} |u(x)| \Big)^p.
\end{equation}
To estimate the second summand on the right-hand side, we note that Caccioppoli's estimate \eqref{caccioppoli} implies
\[
 \sum_{e:|e|\le \frac{R}{2}} | \nabla u(e) |^2 \le C(d,\lambda) R^{-2} \sum_{y:|y|\le R} | u(y) |^2 \le C(d,\lambda)
R^{d-2} \Big( \sup_{x:|x|\le R} |u(x)| \Big)^2.
\]
Together with Jensen's inequality (here we need $q\le2$, that is $p\ge 2$, which is obvious since even $p > d$ from \eqref{C2_sobolev_index}), we obtain that
\begin{align}\nonumber
 \bigg( \sum_{e:|e|\le \frac{R}{2}} | \nabla u(e) |^q \bigg)^{p-1} &\le C(d,p) R^{d(\frac{p}{2}-1)} \bigg( \sum_{e:|e|\le \frac{R}{2}} |
\nabla u(e) |^2 \bigg)^{\frac{p}{2}}\\
 &\le C(d,\lambda,p) R^{d(p-1)-p} \Big( \sup_{x:|x|\le R} |u(x)| \Big)^p.
 \label{C2_W^1q_est}
\end{align}
Substituting \eqref{C2_Lq_est} and \eqref{C2_W^1q_est} into \eqref{C2_grad_estimate2} yields the claim of this step.

\medskip

{\bf Step 3}. Using \eqref{morrey} and bounding the Green function, we conclude that
\begin{equation}\label{C2_step3}
 \Bigg\langle \Bigg( \sup_{u} \frac{\sup_{x:|x|\le R} \frac{|u(x) - u(0)|}{|x|^\alpha}}{\frac{1}{R^\alpha}
\sup_{x:|x|\le R} |u(x)|} \Bigg)^p \Bigg\rangle \le C(d,\lambda,\rho,p,\alpha)
\end{equation}
for all $\alpha < 1$, $p < \infty$, and $R < \infty$, where the outer supremum is taken over all solutions $u(x)$ to
$\nabla^* a \nabla u = 0$ in $\{x:|x|\le R\}$. 
Indeed, Theorem \ref{T} applied to the result \eqref{C2_de_giorgi_green} of Step 2
yields
\begin{multline*}
 \Bigg\langle \Bigg( \sup_{u} \frac{\sup_{x:|x|\le \frac{R}{8}} \frac{|u(x) - u(0)|}{|x|^\alpha}}{\frac{1}{R^\alpha}
\sup_{x:|x|\le R} |u(x)|} \Bigg)^p \Bigg\rangle\\
 \le C(d,\lambda,\rho,p) R^{\alpha p} \Big( R^{d(p-1)-p} \sum_{e:|e|\le \frac{R}{8}} \sum_{b:\frac{R}{4} \le
|b|\le \frac{R}{2}} (|e-b|+1)^{-pd}\\
 + R^{d(p-1)-2p} \sum_{e:|e|\le \frac{R}{8}} \sum_{x:\frac{R}{4} \le |x|\le \frac{R}{2}} (|e-x|+1)^{p(1-d)} \Big)
\end{multline*}
if $\alpha$ and $p$ are related by \eqref{C2_sobolev_index}.
In the domains of $e$ and $b$, we have $|e-b|+1 \ge |b| - |e| \ge R/8$. Therefore the
first double-sum on the right-hand side is bounded by
\[
 C(d,p) R^{2d-pd}.
\]
Likewise the second double-sum is bounded by
\[
 C(d,p) R^{2d+p(1-d)}.
\]
If \eqref{C2_sobolev_index} holds, we thus conclude that
\[
 \Bigg\langle \Bigg( \sup_{u} \frac{\sup_{x:|x|\le \frac{R}{8}} \frac{|u(x) - u(0)|}{|x|^\alpha}}{\frac{1}{R^\alpha}
\sup_{x:|x|\le R} |u(x)|} \Bigg)^p \Bigg\rangle \le C(d,\lambda,\rho,p).
\]
In the region $\{ x : \frac{R}{8} \le |x| \le R \}$, it obviously holds
\[
 \frac{|u(x) - u(0)|}{|x|^\alpha} \le 2 \frac{8^\alpha}{R^\alpha} \sup_{x:|x|\le R} |u(x)|.
\]
Thus we have obtained \eqref{C2_step3} for $p$ and $\alpha$ such that \eqref{C2_sobolev_index} holds. Since in
\eqref{C2_sobolev_index}, $\alpha \to 1$ as $p \to \infty$ and since we can always decrease $p$ and $\alpha$ in the
conclusion \eqref{C2_step3} (in $p$ this follows from Jensen's inequality), the estimate \eqref{C2_step3} indeed holds
for arbitrary $p < \infty$ and $\alpha < 1$.



\subsection{Proof of Lemmas \ref{L.LSI_ss} and \ref{L.LSI_tensor}}\label{ssec:L.LSI}

{\sc Proof of Lemma \ref{L.LSI_ss}}

Without loss of generality, we may assume $\langle \zeta^2 \rangle = 1$. The elementary inequality $\zeta^2 \log \zeta^2 - \zeta^2
+ 1 \le (\zeta^2-1)^2$ then yields
\[
 \langle \zeta^2 \log \zeta^2 \rangle = \langle \zeta^2 \log \zeta^2 - \zeta^2 + 1\rangle \le \langle (\zeta^2-1)^2 \rangle.
\]
Since $(\zeta^2-1)^2 = (|\zeta|-1)^2 (|\zeta|+1)^2$, we find that
\[
 \langle \zeta^2 \log \zeta^2 \rangle \le \langle (|\zeta|+1)^2 \rangle \sup_a (|\zeta|-1)^2 .
\]
Since $\langle \zeta^2 \rangle = 1$, there exists $a_*\in[\lambda,1]$ such that $|\zeta(a_*)|\le 1$. It follows that
\[
  |\zeta(a)|-1 \le |\zeta(a)|-|\zeta(a_*)| \le |\zeta(a) - \zeta(a_*)| \le \osc_a \zeta(a).
\]
Likewise there exists $a^*\in[\lambda,1]$ such that $|\zeta(a^*)|\ge 1$ and therefore
\[
  1-|\zeta(a)| \le |\zeta(a^*)|-|\zeta(a)| \le |\zeta(a^*) - \zeta(a)| \le \osc_a \zeta(a).
\]
Hence it follows that
\[
 \langle \zeta^2\log \zeta^2\rangle \le \langle (|\zeta|+1)^2 \rangle \Big(\osc_a \zeta\Big)^2.
\]
Finally we have that
\[
 \langle (|\zeta|+1)^2 \rangle \le \langle 2\zeta^2+2 \rangle = 4,
\]
and the combination of the previous two inequalities yields \eqref{LSI_ss} with constant $\rho = \frac{1}{8}$.

\bigskip


{\sc Proof of Lemma \ref{L.LSI_tensor}}

The following is a simple adaptation of the usual tensorization proof, cf.\ \cite[Theorem 4.4]{GuionnetZegarlinski}.
Take any enumeration $(e_n)_{n\ge1}$ of the edge set $\E^d$ and denote by $\langle \cdot \rangle_n$ the
$e_n$-marginal of the (product) ensemble $\langle \cdot \rangle$.
We assume that every marginal satisfies the LSI
\[
 \Big\langle \zeta^2 \log \frac{\zeta^2}{\langle \zeta^2 \rangle_n} \Big\rangle_n \le \frac{2}{\rho} \Big( \osc_{a\in[\lambda,1]} \zeta
\Big)^2
\]
for all $\zeta : [\lambda,1] \to \R$. Replacing $\zeta^2$ by $f$ in the definition of the LSI, it suffices to prove
\[
 \Big\langle f \log \frac{f}{\langle f \rangle} \Big\rangle \le \frac{2}{\rho} \sum_{n=1}^\infty \Big\langle \Big( \osc_{a(e_n)}
\sqrt{f} \Big)^2 \Big\rangle
\]
for all positive random variables $f:\Omega\to(0,\infty)$. By a simple density argument, it suffices to consider \emph{local} random
variables, i.e.\ $f$ that depend on $a$ only through a finite number of sites so that the above sum is finite.
We denote iteratively $f_0 := f$ and $f_n := \langle f_{n-1} \rangle_n$. Thus $f_n$ is the average of $f$ over the first $n$
edges. Then the l.-h.\ s.\ of \eqref{LSI} can be expressed as a telescope sum (a finite sum for local random variables):
\begin{align}\nonumber
 \langle f \log f \rangle - \langle f \rangle \log \langle f \rangle &= \sum_{n=1}^\infty \big\langle f_{n-1} \log
f_{n-1} - f_n \log f_n \big\rangle\\ &= \sum_{n=1}^\infty \big\langle \langle f_{n-1} \log
f_{n-1} \rangle_n - \langle f_{n-1} \rangle_n \log \langle f_{n-1} \rangle_n \big\rangle.
\label{LSI_telescope}
\end{align}
The assumption of single-edge LSI yields
\begin{equation}\label{LSI_n}
 \langle f_{n-1} \log f_{n-1} \rangle_n - \langle f_{n-1} \rangle_n \log \langle f_{n-1} \rangle_n \le \frac{2}{\rho} \Big( \osc_{a(e_n)}
\sqrt{f_{n-1}} \Big)^2.
\end{equation}
Notice that the definition of $f_{n-1}$ immediately yields $f_{n-1} = \langle f \rangle_{<n}$, where we have abbreviated the ensemble average over the first $n-1$ edges as $\langle \cdot \rangle_{<n}$. 
We clearly have
\[
 \osc_{a(e_n)} \sqrt{f_{n-1}} = \Big(\sup_{a(e_n)} \langle f \rangle_{<n} \Big)^\frac{1}{2} - \Big(\inf_{a(e_n)} \langle f  \rangle_{<n}\Big)^\frac{1}{2} \le \Big\langle \sup_{a(e_n)} f \Big\rangle_{<n}^\frac{1}{2} -\Big\langle \inf_{a(e_n)} f \rangle_{<n}\Big)^\frac{1}{2}.
\]
By monotonicity of the square root, it follows
\[
 \osc_{a(e_n)} \sqrt{f_{n-1}} \le \Big\langle \Big(\sup_{a(e_n)} \sqrt{f}\Big)^2 \Big\rangle_{<n}^\frac{1}{2} - \Big\langle \Big(\inf_{a(e_n)} \sqrt{f} \Big)^2 \Big\rangle_{<n}^\frac{1}{2}.
\]
Consequently the triangle inequality w.~r.~t.\ $\langle (\cdot)^2 \rangle_{<n}^{\frac{1}{2}}$ on the right-hand side yields
\[
 \osc_{a(e_n)} \sqrt{f_{n-1}} \le \Big\langle \Big(\sup_{a(e_n)} \sqrt{f} - \inf_{a(e_n)} \sqrt{f} \Big)^2\Big\rangle_{<n}^{\frac{1}{2}},
\]
which by definition of $\osc_{a(e_n)}$ can be written as
\begin{equation}
 \osc_{a(e_n)} \sqrt{f_{n-1}} \le \Big\langle \Big( \osc_{a(e_n)} \sqrt{f} \Big)^2 \Big\rangle_{<n}^{\frac{1}{2}}. \label{LSI_n2}
\end{equation}
Finally we collect \eqref{LSI_telescope}, \eqref{LSI_n}, and \eqref{LSI_n2} to obtain
\[
  \langle f \log f \rangle - \langle f \rangle \log \langle f \rangle \le \frac{1}{2\rho} \sum_{n=1}^\infty \Big\langle
\Big( \osc_{a(e_n)} \sqrt{f} \Big)^2 \Big\rangle,
\]
which is the LSI \eqref{LSI} for $f = \zeta^2$.



\begin{thebibliography}{99}

\bibitem{Iwaniec} K.\ Astala, T.\ Iwaniec, G.\ Martin, Elliptic partial differential equations and quasiconformal 
mappings in the plane, {\it Princeton Mathematical Series 48},  Princeton University Press, Princeton, NJ (2009).

\bibitem{Benjaminietal} I.\ Benjamini, H.\ Duminil-Copin, G.\ Kozma and A.\ Yadin, Disorder, entropy and harmonic
functions, Preprint {\tt arXiv:1111.4853} (2011)

\bibitem{Biskup} M.\ Biskup, M.\ Salvi and T.\ Wolff, A central limit theorem for the effective conductance: I. Linear
boundary data and small ellipticity contrasts, Preprint {\tt arXiv:1210.2371} (2012)

\bibitem{CarlenKusuokaStroock} E.\ A.\ Carlen, S.\ Kusuoka, D.\ W.\ Stroock, Upper bounds for symmetric 
Markov transition functions, {\it Ann. Inst. H. Poincar\'e Probab. Statist.} \textbf{23 (2)}, 245-287 (1987)

\bibitem{ConlonNaddaf} J.\ C.\ Conlon, A.\ Naddaf, On homogenization of elliptic
equations with random coefficients,
{\it Electron.\ J.\ Probab.} {\bf 9 (5)}, 1-58 (2000)

\bibitem{DeGiorgi} E.\ De Giorgi, Sulla differenziabilit\`a e l'analiticit\`a delle estremali degli integrali 
multipli regolari, {\it Mem. Accad. Sci. Torino. Cl. Sci. Fis. Mat. Nat. (3)} {\bf 3}, 25-43 (1957)

\bibitem{Delmotte} T.\ Delmotte, In\'egalit\'e de Harnack elliptique sur les graphes, {\it Colloq. Math.} \textbf{72 (1)}, 19--37 (1997)

\bibitem{DeuschelDelmotte} T. Delmotte, J.-D.\ Deuschel, On estimating the derivatives of
symmetric diffusions in stationary random environments, with applications to the $\nabla\phi$ interface model,
{\it Probab. Theory Relat. Fields} {\bf 133}, 358-390 (2005)

\bibitem{Dolzmannetal} G.\ Dolzmann, N.\ Hungerb\"uhler, S.\ M\"uller, Uniqueness and maximal regularity for nonlinear elliptic systems of n-Laplace type with measure
valued right hand side. {\it J.\ Reine Angew.\ Math.} {\bf 520}, 1-35 (2000)

\bibitem{Federbush} P.\ Federbush, Partially alternate derivation of a result by Nelson,
{\it J.\ Math.\ Phys.} {\bf 10 (1)}, 50-52 (1969)

\bibitem{Gloria} A.\ Gloria, Fluctuation of Solutions to Linear Elliptic Equations with Noisy Diffusion Coefficients, {\it Comm. PDE} {\bf 38 (2)}, 304-338 (2012)

\bibitem{GloriaOtto1} A.\ Gloria, F.\ Otto,
An optimal variance estimate in stochastic homogenization of discrete elliptic equations,
{\it Ann.\ Probab.} {\bf 39 (3)}, 779-856 (2011)

\bibitem{GloriaNeukammOtto} A.\ Gloria, S.\ Neukamm, F.\ Otto, Quantification of ergodicity in stochastic homogenization: optimal bounds via
spectral gap on Glauber dynamics, Max Planck Institute for Mathematics in the Sciences Preprint 3/2013

\bibitem{GloriaNeukammOtto2} A.\ Gloria, S.\ Neukamm, F.\ Otto, An optimal quantitative two-scale expansion in stochastic homogenization of
discrete elliptic equations, Max Planck Institute for Mathematics in the Sciences Preprint 41/2013


\bibitem{Gross} L.\ Gross, Logarithmic Sobolev inequalities, {\it American J.\ Math.} {\bf 97 (4)},
1061-1083 (1975)

\bibitem{GuionnetZegarlinski} A.\ Guionnet, B.\ Zegarlinski, Lecture notes on Logarithmic Sobolev
Inequalities, {\it Lecture Notes Math.} {\bf 1801}, 1-134 (2003)

\bibitem{LawlerLimic} G.\ F.\ Lawler and V.\ Limic, Random walk: a modern introduction, {\it Cambridge Stud.\ in
Adv.\ Math.} {\bf 123}, CUP, Cambridge, UK (2010)

\bibitem{Ledoux} M.\ Ledoux, The concentration of measure phenomenon, {\it Math.\ Surveys and Monographs} {\bf 89}, AMS, Providence, RI (2001)

\bibitem{LiebLoss} E.\ H.\ Lieb and M.\ Loss, Analysis, {\it Graduate Stud.\ in Math.} {\bf 14}, 2nd ed., AMS,
Providence, RI (2001)

\bibitem{NaddafSpencerunpublished} A.\ Naddaf, T.\ Spencer,
Estimates on the variance of some homogenization
problems, unpublished

\bibitem{NaddafSpencer} A.\ Naddaf, T.\ Spencer,
On homogenization and scaling limit of some gradient perturbation of a massless free field,
{\it Commun.\ Math.\ Phys.} {\bf 183}, 55-84 (1997)

\bibitem{Nash} J.\ Nash, Continuity of solutions of parabolic and elliptic equations, 
{\it American J.\ Math.} {\bf 80}, 931-954 (1958)

\bibitem{Nelson1} F.\ Nelson, A quartic interaction in two dimensions, in {\it Mathematical theory
of elementary particles} (edited by R.\ Goodman, I.\ Segal), M.\ I.\ T.\ Press (Cambridge, MA), 
69-73 (1966)

\bibitem{Nelson2} F.\ Nelson, The free Markoff field, {\it J.\ Funct.\ Anal.} {\bf 12}, 211-227 (1973)

\bibitem{Nolen} J.\ Nolen, Normal approximation for a random elliptic equation, Preprint (2011). [Available online at
{\tt http://math.duke.edu/{\textasciitilde}nolen/preprints/ellipfluctper\_rev.pdf}.]

\bibitem{RieszNagy} F.\ Riesz and B.\ Sz.-Nagy, Functional Analysis, {\it Dover Books on Adv.\ Math.}, Dover
Publ.\ Inc., New York (1990)

\bibitem{Rossignol} R.\ Rossignol, Noise-stability and central limit theorems for effective resistance of random
electric networks, Preprint, {\tt arXiv:1206.3856} (2012)

\bibitem{StroockZegarlinski} D.\ Stroock, B.\ Zegarlinski, The logarithmic Sobolev inequality
for discrete spin systems on a lattice, {\it Commun.\ Math.\ Phys.} {\bf 149}, 175-193 (1992)

\end{thebibliography}
\end{document}